\documentclass[11pt]{article}
\usepackage{amscd}
\usepackage{amsfonts}
\usepackage{amsmath}
\usepackage{amssymb}
\usepackage{amsthm}
\usepackage{bbm}
\usepackage{CJK}
\usepackage{fancyhdr}
\usepackage{graphicx}
\usepackage{hyperref}
\usepackage{indentfirst}
\usepackage{latexsym}
\usepackage{mathrsfs}
\usepackage{xypic}
\usepackage[compress]{cite}

\newtheorem{theorem}{Theorem}[section]
\newtheorem{lemma}[theorem]{Lemma}
\newtheorem{definition}[theorem]{Definition}
\newtheorem{proposition}[theorem]{Proposition}
\newtheorem{example}[theorem]{Example}

\newtheorem{remark}[theorem]{Remark}

\usepackage[top=1in,bottom=1in,left=1.25in,right=1.25in]{geometry}
\def\<{\langle}
\def\>{\rangle}

\def\o{\otimes}

\begin{document}
\renewcommand{\baselinestretch}{1.2}
\renewcommand{\arraystretch}{1.0}
\title{\bf Modified Rota-Baxter Lie-Yamaguti algebras}

\author{{\bf Wen Teng, Shuangjian Guo\footnote
        { Corresponding author:~~shuangjianguo@126.com}  }\\
{\small  School of Mathematics and Statistics, Guizhou University of Finance and Economics} \\
{\small  Guiyang  550025, P. R. of China}}
 \maketitle
\begin{center}
\begin{minipage}{13.cm}

{\bf \begin{center} ABSTRACT \end{center}}
In this paper, first we introduce the concept of   modified Rota-Baxter Lie-Yamaguti algebras. Then the cohomology of a modified Rota-Baxter Lie-Yamaguti algebra with coefficients in a suitable representation is established.  As   applications, the formal deformations and  abelian extensions of   modified Rota-Baxter Lie-Yamaguti algebras are studied using the second cohomology group.
 \smallskip

 \smallskip

 \end{minipage}
 \end{center}
 \normalsize\vskip0.5cm

\section{Introduction}
\def\theequation{\arabic{section}. \arabic{equation}}
\setcounter{equation} {0}
The concept of Lie-Yamaguti algebras is a generalization of Lie triple systems and Lie algebras, and its root lies in Jacobson's pioneering work \cite{J49,J51}, who formally introduced the concept of Lie triple systems.
Nomizu \cite{N54} further extended this theory and established  invariant affine connections on homogeneous spaces in 1950s.
Building on these foundational ideas, Yamaguti \cite{Y57} introduced the concept of general Lie triple systems, which is a vector space with a pair of bilinear and trilinear operations and  satisfying some intriguing relations to represent the   torsion and curvature tensors of Nomizu's canonical connection. Subsequently, in  \cite{Y67}, Yamaguti proposed the representation and cohomology theory for general Lie triple systems.
Later on, general Lie triple systems were also called   Lie triple algebras by Kikkawa in \cite{K75,K79,K81}.
Kinyon and Weinstein \cite{K01} observed Lie triple algebras, which they called Lie-Yamaguti algebras in their study of Courant algebroids and can be constructed by Leibniz algebras.
Further research on Lie-Yamaguti algebras could be found in \cite{B05,B09, L15, L23,G24,S21,S22,Z22,T23,Z15} and references cited therein.

Inspired by classical $r$-matrix and Poisson structure, Kupershmidt introduced the notion of relative Rota-Baxter operators (also called $\mathcal{O}$-operators) on Lie algebras \cite{K99},  which generalized associative Rota-Baxter operators introduced by Baxter \cite{B60} in the study of fluctuation theory in probability, Rota \cite{R69} are further studied Rota-Baxter operators of  algebra and combination.
Rota-Baxter operators have been widely studied in  many areas of mathematics and physics, including  combinatorics, number theory, operads and quantum field theory \cite{C00}.
Recently, the cohomology and deformation theory of Rota-Baxter operators on various types of algebras have been broken through.
Tang et al. \cite{T19}  developed the  cohomology and deformation of $\mathcal{O}$-operators  on  Lie algebras.  Lazarev  et al.  studied the deformation and cohomology theory of relative Rota-Baxter Lie algebras and its application in triangular Lie bialgebras in \cite{L21}. Later,  Zhao and Qiao discussed the cohomologies and deformations of relative Rota-Baxter operators on Lie-Yamaguti algebras in \cite{Z22}.

Recently, Jiang and Sheng have been considered   cohomologies and deformations of modified $r$-matrices motivated by the deformations of classical $r$-matrices in \cite{J22}.
 Inspired by \cite{J22}, the authors independently studied  cohomology and deformation theory of modified Rota-Baxter Leibniz algebras in  \cite{L22,M22}.  Later, we  studied  cohomology and deformation theory of modified Rota-Baxter Lie triple systems in \cite{GT24}.  It is well known that  Lie-Yamaguti
algebras are a generalization of  Lie algebras, Lie triple systems and  Leibniz algebras, it is very
natural to investgate  modified Rota-Baxter Lie-Yamaguti algebras.  This was our motivation for writing the present paper.  In precisely,  we study the  representation and cohomology   of  modified Rota-Baxter Lie-Yamaguti algebras  and applied them to  the formal deformation  and abelian extension  of  modified Rota-Baxter Lie-Yamaguti algebras.

  This paper is organized as follows.
In Section  \ref{sec:Preliminaries}, we   recall some basic definitions  about Lie-Yamaguti  algebras   and  their cohomology.
 Section \ref{sec:MRBLY}  proposes  the notion of     modified Rota-Baxter Lie-Yamaguti algebras.
 Section \ref{sec:Representations} introduces the representation of modified Rota-Baxter Lie-Yamaguti algebras.
In Section \ref{sec:Cohomology},  we construct the cohomology of the modified Rota-Baxter Lie-Yamaguti algebra.
In Section \ref{sec:deformations},  we use the cohomological approach to study formal deformations of modified Rota-Baxter Lie-Yamaguti algebras.
 Finally, abelian extensions of modified Rota-Baxter Lie-Yamaguti algebras are discussed in Section  \ref{sec:extensions}.

\section{Preliminaries}\label{sec:Preliminaries}
\def\theequation{\arabic{section}. \arabic{equation}}
\setcounter{equation} {0}

Throughout this paper, we work on an algebraically closed field $\mathbb{K}$ of characteristic zero.  In this section, we  recall some basic  definitions of  Lie-Yamaguti  algebra from \cite{K01} and \cite{Y67}.

\begin{definition}
A Lie-Yamaguti  algebra  is a  triple  $(\mathfrak{L}, [-, -], \{-, -, -\})$ in which $\mathfrak{L}$ is a vector space together with  a binary bracket $[-, -]$ and  a ternary bracket $\{-, -, -\}$
on $\mathfrak{L}$   satisfying
\begin{eqnarray*}
&&(LY1)~~ [x, y]=- [y,  x],\\
&&(LY2) ~~\{x,y,z\}=- \{y,x,z\},\\
&&(LY3) ~~ \circlearrowleft_{x,y,z}[[x, y], z]+\circlearrowleft_{x,y,z}\{x, y, z\}=0,\\
&&(LY4) ~~\{[x, y], z, a\}+ \{[z,x],y,a\}+ \{[y,z],x,a\}=0,\\
&&(LY5) ~~\{a, b, [x, y]\}=[\{a, b, x\}, y]+[x,\{a, b, y\}],\\
&&(LY6) ~~ \{a, b, \{x, y, z\}\}=\{\{a, b, x\}, y, z\}+ \{x,  \{a, b, y\}, z\}+ \{x,  y, \{a, b, z\}\},
\end{eqnarray*}
for any $ x, y, z, a, b\in \mathfrak{L}$.
\end{definition}

\begin{example} \label{exam:Lie algebra}
 Let $(\mathfrak{L}, [-, -])$  be a Lie algebra.  Define a  ternary bracket on $\mathfrak{L}$ by
 $$\{x,y,z\}=[[x,y],z], ~~~~\forall x,y,z\in \mathfrak{L}.$$
 Then, $(\mathfrak{L}, [-, -], \{-, -, -\})$ is a   Lie-Yamaguti  algebra.
\end{example}

\begin{example}
 Let $(\mathfrak{L}, [-, -])$  be a Lie algebra  with a reductive decomposition $\mathfrak{L}=\mathfrak{N}\oplus\mathfrak{M}$, i.e.  $[\mathfrak{N},\mathfrak{N}]\subseteq \mathfrak{N}$ and $[\mathfrak{N},\mathfrak{M}]\subseteq \mathfrak{M}$.
Define bilinear bracket $[-,-]_\mathfrak{M}$ and trilinear bracket $\{-,-,-\}_\mathfrak{M}$ on $\mathfrak{M}$ by
the projections of the Lie bracket:
$$[x,y]_\mathfrak{M}=\pi_\mathfrak{M}([x,y]),\{x,y,z\}_\mathfrak{M}=[\pi_\mathfrak{N}([x,y]),z], \forall x,y,z\in \mathfrak{M},$$
where $\pi_\mathfrak{N}:\mathfrak{L}\rightarrow \mathfrak{N},\pi_\mathfrak{M}:\mathfrak{L}\rightarrow \mathfrak{M}$ are the projection map.  Then, $(\mathfrak{M}, [-,-]_\mathfrak{M}, \{-, -, -\}_\mathfrak{M})$ is a   Lie-Yamaguti  algebra.
\end{example}

\begin{example}  \label{exam:Leibniz algebra}
 Let $(\mathfrak{L}, \star)$  be a Leibniz algebra.   Define a binary and ternary bracket on $\mathfrak{L}$ by
 $$[x,y]=x\star y-y\star x, \{x,y,z\}=-(x\star y)\star z, \forall x,y,z\in \mathfrak{L}.$$
 Then, $(\mathfrak{L}, [-, -], \{-, -, -\})$ is a   Lie-Yamaguti  algebra.
\end{example}

\begin{example} \label{exam:2-dimensional LY algebra}
Let $\mathfrak{L}$ be a 2-dimensional  vector space  with a basis $\varepsilon_1,\varepsilon_2$.  If we define  a binary non-zero bracket $[-,-]$ and a ternary non-zero bracket $\{-, -, -\}$ on $\mathfrak{L}$ as follows:
$$[\varepsilon_1,\varepsilon_2]=-[\varepsilon_2,\varepsilon_1]=\varepsilon_1, \{\varepsilon_1,\varepsilon_2, \varepsilon_2\}=-\{\varepsilon_2,\varepsilon_1, \varepsilon_2\}=\varepsilon_1,$$
then $(\mathfrak{L}, [-, -], \{-, -, -\})$ is a   Lie-Yamaguti  algebra.
\end{example}

\begin{example} \label{exam:3-dimensional LY algebra}
Let $\mathfrak{L}$ be a 3-dimensional  vector space  with a basis $\varepsilon_1$, $\varepsilon_2, \varepsilon_3$.  If we define  a binary non-zero bracket $[-,-]$ and a ternary non-zero bracket $\{-, -, -\}$ on $\mathfrak{L}$ as follows:
$$[\varepsilon_1,\varepsilon_2]=-[\varepsilon_2,\varepsilon_1]=\varepsilon_3, \{\varepsilon_1,\varepsilon_2, \varepsilon_1\}=-\{\varepsilon_2,\varepsilon_1, \varepsilon_1\}=\varepsilon_3,$$
then $(\mathfrak{L}, [-, -], \{-, -, -\})$ is a   Lie-Yamaguti  algebra.
\end{example}

\begin{definition}
Let  $(\mathfrak{L}, [-, -], \{-, -, -\})$ be a Lie-Yamaguti  algebra, a linear map $T:\mathfrak{L}\rightarrow \mathfrak{L}$ is said to  a  Rota-Baxter operator of weight -1  if $T$ satisfies
\begin{align}
 [Tx, Ty]=&T([Tx,y]+[x,Ty]-[x,y]),\label{2.1}\\
\{Tx,Ty,Tz\}=&T\big(\{x,Ty,Tz\}+\{Tx,y,Tz\}+\{Tx,Ty,z\}\nonumber\\
&-\{x,y,Tz\}-\{Tx,y,z\}-\{x,Ty,z\}+\{x,y,z\}\big),\label{2.2}
\end{align}
for any $ x, y, z\in \mathfrak{L}$.
\end{definition}

\begin{definition}
 A  Rota-Baxter Lie-Yamaguti algebra of weight -1 is a quadruple $(\mathfrak{L}, [-, -], $ $\{-, -, -\},T)$ consisting of a Lie-Yamaguti  algebra
$(\mathfrak{L}, [-, -], \{-, -, -\})$ and a  Rota-Baxter operator $T$ of weight -1.
\end{definition}

\begin{definition}
Let $(\mathfrak{L}, [-, -], \{-, -, -\})$  be a  Lie-Yamaguti   algebra   and $V$ be a vector space. A representation of $(\mathfrak{L}, [-, -], \{-, -, -\})$ on $V$ consists of a linear map  $\rho: \mathfrak{L}\rightarrow \mathfrak{gl}(V)$ and two bilinear maps $D,\theta:\mathfrak{ L}\times \mathfrak{L}\rightarrow \mathfrak{gl}(V)$ such that
\begin{eqnarray*}
&&(R1)~D(x,y)-\theta(y,x)+\theta(x,y)+\rho([x,y])-\rho(x)\rho(y)+\rho(y)\rho(x)=0,\\
&&(R2)~D([x,y],z)+D([y,z],x)+D([z,x],y)=0,\\
&&(R3)~\theta([x,y],a) =\theta(x,a)\rho(y)-\theta(y,a)\rho(x),\\
&&(R4)~D(a,b)\rho(x)=\rho(x)D(a,b)+\rho(\{a,b,x\}),\\
&&(R5)~\theta(x,[a,b])=\rho(a)\theta(x,b)-\rho(b)\theta(x,a),\\
&&(R6)~D(a,b)\theta(x,y)=\theta(x,y)D(a,b)+\theta(\{a,b,x\},y)+\theta(x,\{a,b,y\}),\\
&&(R7)~\theta(a,\{x,y,z\}) =\theta(y,z)\theta(a,x)-\theta(x,z)\theta(a,y)+D(x,y)\theta(a,z),
\end{eqnarray*}
for all  $x,y,z,a,b\in \mathfrak{L}$. In this case, we also call $V$ a $\mathfrak{L}$-module. We denote a representation by $(V; \rho, \theta, D)$.
\end{definition}
It can be concluded from (R6) that
\begin{eqnarray*}
&&(R6)'~D(a,b)D(x,y)=D(x,y)D(a,b)+D(\{a,b,x\},y)+D(x,\{a,b,y\}).
\end{eqnarray*}

\begin{example}
Let $(\mathfrak{L}, [-, -], \{-, -, -\})$  be a  Lie-Yamaguti   algebra.  We   define linear maps $\mathrm{ad}: \mathfrak{L}\rightarrow \mathfrak{gl}(\mathfrak{L}),\mathcal{L},\mathcal{R}:\otimes^2 \mathfrak{L}\rightarrow \mathfrak{gl}(\mathfrak{L})$ by
\begin{eqnarray*}
& \mathrm{ad}(x)(z):=[x,z],\mathcal{L}(x,y)(z):=\{x,y,z\},\mathcal{R}(x,y)(z):=\{z,x,y\},
\end{eqnarray*}for all $x,y,z\in \mathfrak{L}$.  Then $(\mathfrak{L};\mathrm{ad},\mathcal{L},\mathcal{R})$ forms a representation of $\mathfrak{L}$ on itself,  called the adjoint representation.
\end{example}

Let us recall the Yamaguti cohomology theory on Lie-Yamaguti algebras in \cite{Y67}.  Let $(V; \rho, \theta, D)$ be a
representation of a Lie-Yamaguti algebra $(\mathfrak{L}, [-, -], \{-, -, -\})$, and we denote the set of $(n+1)$-cochains by
$\mathcal{C}^{n+1}_{\mathrm{LY}}(\mathfrak{L},V)$, where
\begin{equation*}
\mathcal{C}_{\mathrm{LY}}^{n+1}(\mathfrak{L},V)= \left\{ \begin{array}{ll}
\mathrm{Hom}(\underbrace{\wedge^2 \mathfrak{L}\otimes\cdots\otimes\wedge^2 \mathfrak{L}}_n,V)\times \mathrm{Hom}(\underbrace{\wedge^2 \mathfrak{L}\otimes\cdots\otimes\wedge^2 \mathfrak{L}}_n\otimes \mathfrak{L},V) &\mbox{ \mbox{}  $n\geq 1,$  }\\
$$\mathrm{Hom}(\mathfrak{L},V)$$ &\mbox{ \mbox{}  $ n=0$.  }
 \end{array}  \right.
\end{equation*}

In the sequel, we recall the coboundary map of $(n+1)$-cochains on a Lie-Yamaguti algebra $\mathfrak{L}$ with the
coefficients in the representation $(V; \rho, \theta, D)$:

If $n\geq 1$, for any $(f,g)\in \mathcal{C}_{\mathrm{LY}}^{n+1}(\mathfrak{L},V)$, $ \mathfrak{K}_i=x_i\wedge y_i\in \wedge^2 \mathfrak{L}, (i=1,2,\cdots,n+1), z\in \mathfrak{L}$,  the coboundary map $\delta^{n+1}=(\delta^{n+1}_I,\delta^{n+1}_{II}):\mathcal{C}_{\mathrm{LY}}^{n+1}(\mathfrak{L},V)\rightarrow \mathcal{C}_{\mathrm{LY}}^{n+2}(\mathfrak{L},V), (f,g)\mapsto (\delta^{n+1}_I(f,g),\delta^{n+1}_{II}(f,g))$ is given as follows:
\begin{align*}
&\delta^{n+1}_I(f,g)(\mathfrak{K}_1,\cdots,\mathfrak{K}_{n+1})\\
=&(-1)^n(\rho(x_{n+1})g(\mathfrak{K}_1,\cdots,\mathfrak{K}_{n},y_{n+1})-\rho(y_{n+1})g(\mathfrak{K}_1,\cdots,\mathfrak{K}_{n},x_{n+1})\\
&-g(\mathfrak{K}_1,\cdots,\mathfrak{K}_{n},[x_{n+1},y_{n+1}]))+\sum_{k=1}^{n}(-1)^{k+1}D(\mathfrak{K}_k)f(\mathfrak{K}_1,\cdots,\widehat{\mathfrak{K}_{k}}\cdots,\mathfrak{K}_{n+1})\\
&+\sum_{1\leq k<l\leq n+1}(-1)^k f(\mathfrak{K}_1,\cdots,\widehat{\mathfrak{K}_{k}}\cdots,\{x_k,y_k,x_l\}\wedge y_l+x_l\wedge \{x_k,y_k,y_l\},\cdots,\mathfrak{K}_{n+1}),
\end{align*}
\begin{align*}
&\delta^{n+1}_{II}(f,g)(\mathfrak{K}_1,\cdots,\mathfrak{K}_{n+1},z)\\
=&(-1)^n(\theta(y_{n+1},z)g(\mathfrak{K}_1,\cdots,\mathfrak{K}_{n},x_{n+1})-\theta(x_{n+1},z)g(\mathfrak{K}_1,\cdots,\mathfrak{K}_{n},y_{n+1}))\\
&+\sum_{k=1}^{n+1}(-1)^{k+1}D(\mathfrak{K}_k)g(\mathfrak{K}_1,\cdots,\widehat{\mathfrak{K}_{k}}\cdots,\mathfrak{K}_{n+1},z)\\
&+\sum_{1\leq k<l\leq n+1}(-1)^k g(\mathfrak{K}_1,\cdots,\widehat{\mathfrak{K}_{k}}\cdots,\{x_k,y_k,x_l\}\wedge y_l+x_l\wedge \{x_k,y_k,y_l\},\cdots,\mathfrak{K}_{n+1},z)\\
&+\sum_{k=1}^{n+1}(-1)^kg(\mathfrak{K}_1,\cdots,\widehat{\mathfrak{K}_{k}}\cdots,\mathfrak{K}_{n+1},\{x_k,y_k,z\}),
\end{align*}
where $~\widehat{}~$ denotes omission.
For the case that $n=0$, for any $f\in \mathcal{C}_{\mathrm{LY}}^1(\mathfrak{L},V)$ , the coboundary map
 $\delta^1=(\delta^1_I,\delta^1_{II})\text{:}$\, $\mathcal{C}_{\mathrm{LY}}^1(\mathfrak{L},V)\rightarrow \mathcal{C}_{\mathrm{LY}}^2(\mathfrak{L},V),f\rightarrow (\delta^1_I(f),\delta^1_{II}(f))$ is given by:
\begin{align*}
\delta^1_I(f)(x,y)=&\rho(x)f(y)-\rho(y)f(x)-f([x,y]),\\
\delta^1_{II}(f)(x,y,z)=&D(x,y)f(z)+\theta(y,z)f(x)-\theta(x,z)f(y)-f(\{x,y,z\}).
\end{align*}
The corresponding cohomology groups are denoted by $\mathcal{H}^{\ast}_{\mathrm{LY}}(\mathfrak{L},V).$

\section{Modified Rota-Baxter Lie-Yamaguti algebras}\label{sec:MRBLY}
\def\theequation{\arabic{section}. \arabic{equation}}
\setcounter{equation} {0}

In this section, we introduce the   concept of  modified Rota-Baxter Lie-Yamaguti algebras.
We explore the relationship between modified Rota-Baxter operator and Rota-Baxter operator of weight -1 and Nijenhuis operator respectively,
and provide some   examples.

Motivated by the   modified $r$-matrix\cite{J22}, we propose the definition of a  modified Rota-Baxter operator on a Lie-Yamaguti  algebra.

\begin{definition}
Let  $(\mathfrak{L}, [-, -], \{-, -, -\})$ be a Lie-Yamaguti  algebra, a linear map $R:\mathfrak{L}\rightarrow \mathfrak{L}$ is said to  a  modified Rota-Baxter operator  if $R$ satisfies
\begin{align}
 [Rx, Ry]=&R([Rx,y]+[x,Ry])-[x,y]\label{3.1},\\
\{Rx,Ry,Rz\}=&R(\{x,Ry,Rz\}+\{Rx,y,Rz\}+\{Rx,Ry,z\}+\{x,y,z\})\nonumber\\
&-\{Rx,y,z\}-\{x,Ry,z\}-\{x,y,Rz\},\label{3.2}
\end{align}
for any $ x, y, z\in \mathfrak{L}$.
\end{definition}

\begin{definition}
 A modified Rota-Baxter Lie-Yamaguti algebra is a quadruple $(\mathfrak{L}, [-, -], $ $\{-, -, -\},R)$ consisting of a Lie-Yamaguti  algebra
$(\mathfrak{L}, [-, -], \{-, -, -\})$ and a modified Rota-Baxter operator $R$.
\end{definition}

\begin{definition}
A homomorphism between two  modified Rota-Baxter Lie-Yamaguti algebras $(\mathfrak{L}, [-, -], $ $\{-, -, -\},R)$ and $(\mathfrak{L}', [-, -]',  \{-, -, -\}',R')$ is a linear map $\varphi: \mathfrak{L}\rightarrow \mathfrak{L}'$ satisfying
\begin{eqnarray*}
&&\varphi([x, y])=[\varphi(x), \varphi(y)]', ~\varphi(\{x, y, z\})=\{\varphi(x), \varphi(y),\varphi(z)\}' \   \mathrm{and}\    R'\circ\varphi=\varphi\circ R,
\end{eqnarray*}
 for all $x,y,z\in \mathfrak{L}$.
\end{definition}

\begin{remark}
(i) When a   modified Rota-Baxter  Lie-Yamaguti algebra  reduces to a Modified Rota-Baxter Lie triple system, that is  $[-, -]=0$, we get the notion of a  modified Rota-Baxter   Lie
triple system immediately.
See \cite{GT24} for more details about  modified  Rota-Baxter Lie  triple systems.
\\
(ii) When a Lie-Yamaguti algebra reduces to a Lie algebra, that is  $\{-, -, -\}=0$, we get the notion of a  modified $r$-matric. See \cite{J22} for more details about  modified $r$-matrices.
\end{remark}

\begin{example}
Let $(\mathfrak{L}, [-, -], \{-, -, -\})$ be the 2-dimensional Lie-Yamaguti algebra given in $\mathrm{Example}$ \ref{exam:2-dimensional LY algebra}.
Then, for $k,k_1\in \mathbb{K},$
$$R=\left(
        \begin{array}{cc}
          1 & k_1  \\
          0 & k
        \end{array}
      \right)$$
is a modified Rota-Baxter operator.
\end{example}

\begin{example}
Let $(\mathfrak{L}, [-, -], \{-, -, -\})$ be the 3-dimensional Lie-Yamaguti algebra given in $\mathrm{Example}$ \ref{exam:3-dimensional LY algebra}.
Then, $k,k_1,k_2,k_3\in \mathbb{K},$
$$R=\left(
        \begin{array}{ccc}
          0 & k_1 & 0  \\
          \frac{1-k^2}{k_1} & k  & 0 \\
          k_2 & k_3  & k
        \end{array}
      \right)$$
is a modified Rota-Baxter operator.
\end{example}

\begin{example}
An identity map $\mathrm{id}_\mathfrak{L}:\mathfrak{L}\rightarrow \mathfrak{L}$   is a modified Rota-Baxter operator.
\end{example}

\begin{proposition}
Let  $(\mathfrak{L}, [-, -], \{-, -, -\})$ be a Lie-Yamaguti  algebra. A linear map $R$ is    a  modified Rota-Baxter operator
 if and only if  $-R$ is also   a  modified Rota-Baxter operator.
\end{proposition}

There is a close relationship between  Rota-Baxter operators of weight -1 and  modified Rota-Baxter operators.

\begin{proposition}
Let  $(\mathfrak{L}, [-, -], \{-, -, -\})$ be a Lie-Yamaguti  algebra.  If $T$ is    a    Rota-Baxter operator of weight -1,
then  $2T-\mathrm{id}_\mathfrak{L}$ is    a  modified Rota-Baxter operator.
\end{proposition}

\begin{proof}
For any $x, y, z\in  \mathfrak{L}$, by Eqs. \eqref{2.1} and  \eqref{2.2},  we have
\begin{align*}
&[(2T-\mathrm{id}_\mathfrak{L})x, (2T-\mathrm{id}_\mathfrak{L})y]\\
&\ =[2Tx- x, 2Ty-y]\\
&\ =4[Tx, Ty]-2[Tx, y]-2[x, Ty]+[x, y]\\
&\ =4T([Tx,y]+[x,Ty]-[x,y])-2[Tx, y]-2[x, Ty]+[x, y]\\
&\ =(2T-\mathrm{id}_\mathfrak{L})([(2T-\mathrm{id}_\mathfrak{L})x,y]+[x,(2T-\mathrm{id}_\mathfrak{L})y])-[x,y],\\
& \{(2T-\mathrm{id}_\mathfrak{L})(x),  (2T-\mathrm{id}_\mathfrak{L})(y), (2T- \mathrm{id}_\mathfrak{L})(z)\}\}\\
&\ =\{2Tx- x,  2Ty- y, 2Tz- z\}\\
&\ = 8\{Tx, Ty, Tz\}-4\{Tx, Ty, z]-4[T(x), y, Tz\}-4\{x, Ty, Tz\}\\
&\ \ \ \ +2\{x, y, Tz\}+2\{x, Ty, z\}+2\{Tx, y, z\}-\{x, y, z\}\\
&\ = 8T(\{Tx, Ty, z\} + \{Tx, y, Tz\} + \{x, Ty, Tz\}-\{Tx, y, z\} - \{x, Ty, z\} - \{x, y, Tz\}+\{x, y, z\})\\
&\ \ \ \ -4\{Tx, Ty, z\}-4\{Tx, y, Tz\}-4\{x, Ty, Tz\}\\
&\ \ \ \ +2\{x, y, Tz\}+2\{x, Ty, z\}+2\{Tx, y, z\}-\{x, y, z\}\\
&\ = (2T- \mathrm{id}_\mathfrak{L})( \{(2T- \mathrm{id}_\mathfrak{L})x, (2T- \mathrm{id}_\mathfrak{L})y, z\}+\{x, (2T- \mathrm{id}_\mathfrak{L})y, (2T- \mathrm{id}_\mathfrak{L})z\}\\
&\ \ \ \ +\{(2T- \mathrm{id}_\mathfrak{L})x, y, (2T- \mathrm{id}_\mathfrak{L})z\}\}+\{x, y, z\})- \{(2T-\mathrm{id}_\mathfrak{L})x, y, z\}-\{x, (2T-\mathrm{id}_\mathfrak{L})y, z\}\\
&\ \ \ \ -\{x, y, (2T- \mathrm{id}_\mathfrak{L})z\}.
\end{align*}
The proposition follows.
\end{proof}
Recall from~\cite{S21} that a Nijenhuis operator on a Lie-Yamaguti algebra  $(\mathfrak{L}, [\cdot, \cdot], \{\cdot, \cdot, \cdot\})$ is a
linear map $N:\mathfrak{L}\rightarrow \mathfrak{L}$
 satisfies
\begin{align*}
 [Nx,Ny]=&N([Nx, y]+[x, Ny]-N[x, y]),\\
\{Nx,Ny,Nz\}=&N(\{Nx,Ny,z\}+\{Nx,y,Nz\}+\{x,Ny,Nz\})-\\
&N^{2}(\{Nx,y,z\}+\{x,Ny,z\}+\{x,y,Nz\})+N^{3}\{x,y,z\},
\end{align*}
for all $x,y,z\in \mathfrak{L}.$ The relationship between the modified Rota-Baxter operator and Nijenhuis operator is as follows, which proves to be obvious.

\begin{proposition}
Let  $(\mathfrak{L}, [-, -], \{-, -, -\})$ be a Lie-Yamaguti  algebra and $N$ a Nijenhuis operator.  If $N^2=\mathrm{id}_\mathfrak{L}$, then  $N$ is    a    Nijenhuis operator
 if and only if  $N$ is    a  modified Rota-Baxter operator.
\end{proposition}

\section{ Representations of modified Rota-Baxter Lie-Yamaguti algebras}\label{sec:Representations}
\def\theequation{\arabic{section}. \arabic{equation}}
\setcounter{equation} {0}

In this section, we introduce the  representation of  modified Rota-Baxter Lie-Yamaguti algebras.
Furthermore, we discuss its relationship with the representation of Rota-Baxter Lie-Yamaguti algebras of weight -1.
 Finally, we establish a new modified Rota-Baxter Lie-Yamaguti algebra $\mathfrak{L}_R$ and give its representations, which prepares for introducing the cohomology of the modified Rota-Baxter Lie-Yamaguti algebra in the next section.

\begin{definition}
  A representation of the modified Rota-Baxter  Lie-Yamaguti algebra $(\mathfrak{L}, $ $[-, -],  \{-, -, -\},R)$  is a quintuple $(V; \rho,   \theta, D, R_V)$ such that the following conditions are satisfied:

  (i) $(V; \rho,   \theta,D)$ is a representation of the  Lie-Yamaguti algebra $(\mathfrak{L}, [-, -],  \{-, -, -\})$;

  (ii)  $R_V:V\rightarrow V$  is a linear map satisfying the following equations
\begin{align}
\rho(Rx)R_Vu=&R_V(\rho(Rx)u+\rho(x)R_Vu)-\rho(x)u,\label{4.1}\\
\theta(Rx, Ry)R_Vu=&R_V\big(\theta(Rx, Ry) u +\theta(Rx, y)R_Vu+\theta(x, Ry)R_Vu+\theta(x, y)u\big)\nonumber\\
&- \theta(Rx, y)u-\theta(x, y)R_Vu-\theta(x, Ry)u,\label{4.2}
\end{align}
for any $x,y\in \mathfrak{L}$ and $ u\in V.$
\end{definition}

It can be concluded from  Eq. \eqref{4.2}  that
\begin{align}
D(Rx, Ry)R_Vu=&R_V\big(D(Rx, Ry) u +D(Rx, y)R_Vu+D(x, Ry)R_Vu+D(x, y)u\big)\nonumber\\
&- D(Rx, y)u-D(x, y)R_Vu-D(x, Ry)u,\label{4.3}
\end{align}

\begin{example}
$(\mathfrak{L};\mathrm{ad},\mathcal{L},\mathcal{R},R)$ is an adjoint representation of the modified Rota-Baxter  Lie-Yamaguti algebra $(\mathfrak{L},  [-, -],  \{-, -, -\}, R)$.
\end{example}

The representation of  modified  Rota-Baxter  Lie-Yamaguti algebras is closely related to the representation of  Rota-Baxter  Lie-Yamaguti algebras of   weight -1.

\begin{definition}
  A representation of the   Rota-Baxter  Lie-Yamaguti algebra $(\mathfrak{L}, $ $[-, -], $ $ \{-, -, -\}, T)$ of   weight -1 is a quintuple $(V; \rho,   \theta, D, T_V)$ such that the following conditions are satisfied:

  (i) $(V; \rho,   \theta,D)$ is a representation of the  Lie-Yamaguti algebra $(\mathfrak{L}, [-, -],  \{-, -, -\})$;

  (ii)  $T_V:V\rightarrow V$  is a linear map satisfying the following equations
\begin{align}
\rho(Tx)T_V(u)=&T_V\big(\rho(Tx)u+\rho(x)T_V(u)-\rho(x)u\big),\label{4.4}\\
\theta(Tx, Ty)T_Vu=&T_V\big(\theta(Tx, Ty) u +\theta(Tx, y)T_Vu+\theta(x, Ty)T_Vu\nonumber\\
&-\theta(Tx, y) u -\theta(x, Ty)u-\theta(x, y)T_Vu+\theta(x, y)u\big),\label{4.5}
\end{align}
for any $x,y\in \mathfrak{L}$ and $ u\in V.$
\end{definition}

\begin{proposition}
If  $(V; \rho,   \theta, D, T_V)$  is a representation of the  Rota-Baxter  Lie-Yamaguti algebra $(\mathfrak{L}, $ $[-, -], $ $ \{-, -, -\}, T)$  of   weight -1,
then $(V; \rho,   \theta, D, 2T_V-\mathrm{id}_V)$  is a representation of the modified  Rota-Baxter  Lie-Yamaguti algebra $(\mathfrak{L}, $ $[-, -], $ $ \{-, -, -\}, 2T-\mathrm{id}_\mathfrak{L})$,
\end{proposition}

\begin{proof}
 For any $x,y\in \mathfrak{L}$ and $u\in V$, in the light of  Eqs. \eqref{4.4} and  \eqref{4.5}, we have
 \begin{align*}
 &\rho((2T-\mathrm{id}_\mathfrak{L})x)(2T_V-\mathrm{id}_V)(u)\\
 &\  =\rho((2Tx-x))(2T_Vu-u)\\
  &\  =4\rho(Tx)T_Vu-2\rho(Tx)u-2\rho(x)T_Vu+\rho(x)u\\
  &\  =4T_V\big(\rho(Tx)u+\rho(x)T_V(u)-\rho(x)u\big)-2\rho(Tx)u-2\rho(x)T_Vu+\rho(x)u\\
  &\  =( 2T_V-\mathrm{id}_V)(\rho(( 2T-\mathrm{id}_\mathfrak{L})x)u+\rho(x)( 2T_V-\mathrm{id}_V)(u))-\rho(x)u,\\
 & \theta((2T-\mathrm{id}_\mathfrak{L})x, (2T-\mathrm{id}_\mathfrak{L})y)(2T_V-\mathrm{id}_V)(u)\\
&\ = \theta(2Tx- x, 2Ty- y)(2T_V u- u)\\
&\ = 8\theta(Tx, Ty)T_V(u)-4\theta(x, Ty)T_V (u) -4\theta(Tx, y)T_Vu -4\theta(Tx, Ty)u\\
&\ \ \ \ +2\theta(x, y)T_Vu+2\theta(x, Ty)u+2\theta(Tx, y)u -\theta(x, y)u\\
&\ = 8 T_V\big(\theta(Tx, Ty) u +\theta(Tx, y)T_V(u)+\theta(x, Ty)T_V(u)-\theta(Tx, y) u -\theta(x, Ty)u-\theta(x, y)T_Vu \\
&\ \ \ \ +\theta(x, y)u\big) -4\theta(x, Ty)T_V (u) -4\theta(Tx, y)T_Vu -4\theta(Tx, Ty)u+2\theta(x, y)T_Vu+2\theta(x, Ty)u\\
&\ \ \ \ +2\theta(Tx, y)u -\theta(x, y)u\\
&\ =(2T_V-\mathrm{id}_V)\big(\theta((2T-\mathrm{id}_\mathfrak{L})x , (2T-\mathrm{id}_\mathfrak{L})y) u +\theta((2T-\mathrm{id}_\mathfrak{L})x, y)(2T_V-\mathrm{id}_V)u\\
&\ \ \ \ +\theta(x, (2T-\mathrm{id}_\mathfrak{L})y)(2T_V-\mathrm{id}_V)u+\theta(x, y)u\big)- \theta((2T-\mathrm{id}_\mathfrak{L})x, y)u-\theta(x, y)(2T_V-\mathrm{id}_V)u\\
&\ \ \ \ -\theta(x, (2T-\mathrm{id}_\mathfrak{L})y)u.
\end{align*}
The proposition follows.
\end{proof}

Next we  construct the semidirect product of the modified  Rota-Baxter  Lie-Yamaguti algebra.

\begin{proposition}
If  $(V; \rho,   \theta, D, R_V)$  is a representation of the  modified  Rota-Baxter  Lie-Yamaguti algebra $(\mathfrak{L}, $ $[-, -], $ $ \{-, -, -\}, R)$, then $\mathfrak{L} \oplus V$ is a modified  Rota-Baxter  Lie-Yamaguti   algebra  under the following maps:
\begin{eqnarray*}
&&[x+u, y+v]_{\ltimes}:=[x, y]+\rho(x)(v)-  \rho(y)(u),\\
&&\{x+u, y+v, z+w\}_{\ltimes}:=\{x, y, z\}+D(x, y)(w)- \theta(x, z)(v)+ \theta(y, z)(u),\\
&&R\oplus R_V(x+u)=Rx+R_Vu,
\end{eqnarray*}
for all  $x, y, z\in \mathfrak{L}$ and $u, v, w\in V$. In the case, the modified  Rota-Baxter  Lie-Yamaguti   algebra $\mathfrak{L} \oplus V$  is called a semidirect product of $\mathfrak{L}$ and $V$, denoted by $\mathfrak{L}\ltimes V=(\mathfrak{L} \oplus V,[-,-]_{\ltimes},\{-,-,-\}_{\ltimes},R\oplus R_V)$.
\end{proposition}

\begin{proof}
In view of \cite{Y67}, $(\mathfrak{L} \oplus V,[-,-]_{\ltimes},\{-,-,-\}_{\ltimes})$ is a Lie-Yamaguti algebra.

Next, for any $x, y, z\in \mathfrak{L}$ and $u, v, w\in V$,   in view of  Eqs. \eqref{3.1}, \eqref{3.2}  and     \eqref{4.1}-\eqref{4.3},  we have
\begin{align*}
&[R\oplus R_V(x+u), R\oplus R_V(y+v)]_{\ltimes}\\
&\ =[Rx+R_Vu, Ry+ R_Vv]_{\ltimes}\\
&\ =[Rx, Ry]+\rho(Rx)(R_Vv)-  \rho(Ry)(R_Vu)\\
&\ =R([Rx,y]+[x,Ry])-[x,y]+R_V(\rho(Rx)v+\rho(x)R_Vv)-\rho(x)v\\
&\ \ \ \ -  R_V(\rho(Ry)u+\rho(y)R_Vu)+\rho(y)u\\
&\ =(R\oplus R_V)([x+u, (R\oplus R_V)(y+v)]_{\ltimes}+[(R\oplus R_V)(x+u), y+v]_{\ltimes})-[x+u, y+v]_{\ltimes},\\
& [R\oplus R_V(x+u), R\oplus R_V(y+v), R\oplus R_V(z+w)\}_{\ltimes}\\
&\ = \{Rx, Ry, Rz\} +D(Rx, Ry)R_V(w)-\theta(Rx, Rz)R_V(v)+\theta(Ry, Rz)R_V(u)\\
&\ =R( \{Rx, Ry, z\}+\{x, Ry, Rz\}+\{Rx, y, Rz\}+\{x, y, z\})- \{Rx, y, z\}-\{x, Ry, z\}-\{x, y, Rz\}\\
&\ \ \ \  +R_V\big(D(Rx, Ry) w +D(Rx, y)R_V(w)+D(x, Ry)R_V(w)+D(x, y)w\big)\\
&\ \ \ \ - D(Rx, y)w-D(x, y)R_V(w)-D(x, Ry)w\\
&\ \ \ \ -R_V\big(\theta(Rx, Rz) v +\theta(Rx, z)R_V(v)+\theta(x, Rz)R_V(v)+\theta(x, z)v\big)\\
&\ \ \ \ + \theta(Rx, z)v+\theta(x, z)R_V(v)+\theta(x, Rz)v\\
&\ \ \ \  +R_V\big(\theta(Ry, Rz)u +\theta(Ry, z)R_V(u)+\theta(y, Rz)R_V(u)+\theta(y, z)u\big)\nonumber\\
&\ \ \ \ - \theta(Ry, z)u-\theta(y, z)R_V(u)-\theta(y, Rz)u\\
&\ =R\oplus R_V( \{R\oplus R_V(x+u), R\oplus R_V(y+v), z+w\}_{\ltimes}+\{x+u, R\oplus R_V(y+v), R\oplus R_V(z+w)\}_{\ltimes}\\
&\ \ \ \ +\{R\oplus R_V(x+u), y+v, R\oplus R_V(z+w)\}_{\ltimes}+\{x+u, y+v,z+w\}_{\ltimes})\\
&\ \ \ \  - \{R\oplus R_V(x+u), y+v, z+w\}_{\ltimes}-\{x+u, R\oplus R_V(y+v), z+w\}_{\ltimes}\\
&\ \ \ \ -\{x+u, y+v, R\oplus R_V(z+w)\}_{\ltimes}.
\end{align*}
Therefore, $(\mathfrak{L} \oplus V,[-,-]_{\ltimes},\{-,-,-\}_{\ltimes},R\oplus R_V)$ is a modified  Rota-Baxter  Lie-Yamaguti   algebra.
\end{proof}

\begin{proposition}\label{prop:nLY}
Let $(\mathfrak{L}, $ $[-, -], $ $ \{-, -, -\}, R)$  be  a    modified  Rota-Baxter  Lie-Yamaguti algebra, Define   new operations as follows:
\begin{align}
 [x, y]_R=&[Rx,y]+[x,Ry],\label{4.6}\\
\{x,y,z\}_R=&\{x,Ry,Rz\}+\{Rx,y,Rz\}+\{Rx,Ry,z\}+\{x,y,z\}, \forall x,y,z\in \mathfrak{L}.\label{4.7}
\end{align}
 Then, (i) $(\mathfrak{L},  [-, -]_R,  \{-, -, -\}_R)$ is a   Lie-Yamaguti algebra. We denote this Lie-Yamaguti algebra by $\mathfrak{L}_R.$\\
 (ii) $(\mathfrak{L}_R, R)$ is a modified  Rota-Baxter  Lie-Yamaguti algebra.
\end{proposition}

\begin{proof}
Both (i) and (ii) are direct verification, so we have omitted the details.
\end{proof}

\begin{proposition}\label{prop:nLYr}
Let $(V; \rho,   \theta, D, R_V)$  be a representation of the  modified  Rota-Baxter  Lie-Yamaguti algebra $(\mathfrak{L}, $ $[-, -], $ $ \{-, -, -\}, R)$.
Define   linear maps $\rho_R: \mathfrak{L} \rightarrow \mathfrak{gl}(V), \theta_R, D_R: \mathfrak{L}\times \mathfrak{L} \rightarrow \mathfrak{gl}(V)$  by
 \begin{align}
 \rho_R(x)u:=&\rho(Rx)u-R_V(\rho(x)u), \label{4.8}\\
\theta_R(x, y)u:=&\theta(Rx, Ry) u+\theta(Rx, y) u+\theta(x, Ry) u-R_V(\theta(x, y)u), \label{4.9}
\end{align}
for any $x,y\in \mathfrak{L}$ and $u\in V.$
It can be concluded from  Eq. \eqref{4.9}  that
 \begin{align}
D_R(x, y)u=&D(Rx, Ry) u+D(Rx, y) u+D(x, Ry) u-R_V(D(x, y)u). \label{4.10}
\end{align}
Then $(V; \rho_R,   \theta_R, D_R)$  is a representation of  $\mathfrak{L}_R$. Moreover, $(V; \rho_R,   \theta_R, D_R,R_V)$  is a representation of the  modified  Rota-Baxter  Lie-Yamaguti algebra  $(\mathfrak{L}_R,R)$.
\end{proposition}

\begin{proof}
First, through direct verification, $(V; \rho_R,   \theta_R, D_R)$  is a representation of the  Lie-Yamaguti algebra $\mathfrak{L}_R$.
Further, for any $x, y, z\in \mathfrak{L}$ and $u, v, w\in V$,   by  Eqs.   \eqref{4.1}-\eqref{4.3},  we have
\begin{align*}
&\rho_R(Rx)R_Vu\\
&\ =\rho(R^2x)R_Vu-R_V(\rho(Rx)R_Vu)\\
&\ =R_V(\rho(R^2x)u+\rho(Rx)R_Vu)-\rho(Rx)u-R^2_V(\rho(Rx)u+\rho(x)R_Vu)+R_V(\rho(x)u)\\
&\ =R_V\big(\rho(R^2x)u-R_V(\rho(Rx)u)+\rho(Rx)R_Vu- R_V(\rho(x)R_Vu)\big)-(\rho(Rx)u-R_V(\rho(x)u))\\
&\ =R_V(\rho_R(Rx)u+\rho_R(x)R_Vu)-\rho_R(x)u,
\end{align*}
\begin{align*}
& \theta_R(Rx, Ry)R_{V}(v)\\
\quad &= \theta(R^2x, R^2y) R_{V}(v)+\theta(R^2x, Ry) R_{V}(v)+\theta(Rx, R^2y) R_{V}(v)-R_V(\theta(Rx, Ry)R_{V}(v))\\
\quad &= R_V\big(\theta(R^2x, R^2y) v +\theta(R^2x, Ry)R_V(v)+\theta(Rx, R^2y)R_V(v)+\theta(Rx, Ry)v\big)\\
\quad &\ \ \ \ - \theta(R^2x, Ry)v-\theta(Rx, Ry)R_V(v)-\theta(Rx, R^2y)v\\
\quad &\ \ \ \  +R_V\big(\theta(R^2x, Ry) v +\theta(R^2x, y)R_V(v)+\theta(R, Ry)R_V(v)+\theta(Rx, y)v\big)\\
\quad &\ \ \ \  - \theta(R^2x, y)v-\theta(Rx, y)R_V(v)-\theta(Rx, Ry)v\\
\quad &\ \ \ \  +R_V\big(\theta(Rx, R^2y) v +\theta(Rx, Ry)R_V(v)+\theta(x, R^2y)R_V(v)+\theta(x, Ry)v\big)\\
\quad &\ \ \ \  - \theta(Rx, Ry)v-\theta(x, Ry)R_V(v)-\theta(x, R^2y)v\\
\quad &\ \ \ \  -R_V\big(R_V\big(\theta(Rx, Ry) v +\theta(Rx, y)R_V(v)+\theta(x, Ry)R_V(v)+\theta(x, y)v\big)\\
\quad &\ \ \ \  - \theta(Rx, y)v-\theta(x, y)R_V(v)-\theta(x, Ry)v \big)\\
\quad &=R_V\big(\theta_R(Rx, Ry) v +\theta_R(Rx, y)R_V(v)+\theta_R(x, Ry)R_V(v)+\theta_R(x, y)v\big)\\
\quad &\ \ \ \  - \theta_R(Rx, y)v-\theta_R(x, y)R_V(v)-\theta_R(x, Ry)v.
\end{align*}
Hence, $(V; \rho_R,   \theta_R, D_R,R_V)$  is a representation of the  modified  Rota-Baxter  Lie-Yamaguti algebra  $(\mathfrak{L}_R,R)$.
\end{proof}
\begin{example}
$(\mathfrak{L};\mathrm{ad}_R,\mathcal{L}_R,\mathcal{R}_R,R)$ is an adjoint representation of the modified Rota-Baxter  Lie-Yamaguti algebra $(\mathfrak{L}_R, R)$,
where
\begin{align*}
\mathrm{ad}_R(x)(z):=&[Rx,z]-R[x,z],\\
\mathcal{L}_R(x,y)(z):=&\{Rx,Ry,z\}+\{Rx,y,z\}+\{x,Ry,z\}-R\{x,y,z\},\\
\mathcal{R}_R(x,y)(z):=&\{z,Rx,Ry\}+\{z,Rx,y\}+\{z,Rx,y\}-R\{z,x,y\},
\end{align*}for any $x,y,z\in \mathfrak{L}_R$.
\end{example}

\section{ Cohomology of modified Rota-Baxter Lie-Yamaguti algebras}\label{sec:Cohomology}
\def\theequation{\arabic{section}. \arabic{equation}}
\setcounter{equation} {0}

In this section, we will construct the cohomology of  modified Rota-Baxter Lie-Yamaguti algebras.

 Firstly, the cohomology of the modified Rota-Baxter operator is given by  Yamaguti cohomology \cite{Y67}.

Let $(V; \rho,   \theta, D, R_V)$  be a representation of the  modified  Rota-Baxter  Lie-Yamaguti algebra $(\mathfrak{L}, $ $[-, -], $ $ \{-, -, -\}, R)$.
Recall that Proposition \ref{prop:nLY} and Proposition \ref{prop:nLYr} give a new Lie-Yamaguti algebra $\mathfrak{L}_R$ and
a new representation $(V; \rho_R,   \theta_R, D_R)$ over $\mathfrak{L}_R$. Consider the cochain complex of $\mathfrak{L}_R$ with coefficients in $V$:
\begin{equation*}
(\mathcal{C}_{\mathrm{LY}}^{\ast}(\mathfrak{L}_R,V),\partial)=(\oplus_{n=0}^{\infty}\mathcal{C}_{\mathrm{LY}}^{n+1}(\mathfrak{L}_R,V),\partial).
\end{equation*}
More precisely,
\begin{equation*}
\mathcal{C}_{\mathrm{LY}}^{n+1}(\mathfrak{L}_R,V)= \left\{ \begin{array}{ll}
\mathrm{Hom}(\underbrace{\wedge^2 \mathfrak{L}_R\otimes\cdots\otimes\wedge^2 \mathfrak{L}_R}_n,V)\times \mathrm{Hom}(\underbrace{\wedge^2 \mathfrak{L}_R\otimes\cdots\otimes\wedge^2 \mathfrak{L}_R}_n\otimes \mathfrak{L}_R,V) &\mbox{ \mbox{}  $n\geq 1,$  }\\
$$\mathrm{Hom}(\mathfrak{L}_R,V)$$ &\mbox{ \mbox{}  $ n=0$  }
 \end{array}  \right.
\end{equation*}
 and its coboundary map $\partial^{n+1}=(\partial^{n+1}_I,\partial^{n+1}_{II}):\mathcal{C}_{\mathrm{LY}}^{n+1}(\mathfrak{L}_R,V)\rightarrow \mathcal{C}_{\mathrm{LY}}^{n+2}(\mathfrak{L}_R,V), (f,g)\mapsto (\partial^{n+1}_I(f,g),\partial^{n+1}_{II}(f,g))$ is given as follows:
\begin{align*}
&\partial^{n+1}_I(f,g)(\mathcal{K}_1,\cdots,\mathcal{K}_{n+1})\\
=&(-1)^n(\rho_R(x_{n+1})g(\mathcal{K}_1,\cdots,\mathcal{K}_{n},y_{n+1})-\rho_R(y_{n+1})g(\mathcal{K}_1,\cdots,\mathcal{K}_{n},x_{n+1})\\
&-g(\mathcal{K}_1,\cdots,\mathcal{K}_{n},[x_{n+1},y_{n+1}]_R))+\sum_{k=1}^{n}(-1)^{k+1}D_R(\mathcal{K}_k)f(\mathcal{K}_1,\cdots,\widehat{\mathcal{K}_{k}}\cdots,\mathcal{K}_{n+1})\\
&+\sum_{1\leq k<l\leq n+1}(-1)^k f(\mathcal{K}_1,\cdots,\widehat{\mathcal{K}_{k}}\cdots,\{x_k,y_k,x_l\}_R\wedge y_l+x_l\wedge \{x_k,y_k,y_l\}_R,\cdots,\mathcal{K}_{n+1}),
\end{align*}
\begin{align*}
&\partial^{n+1}_{II}(f,g)(\mathcal{K}_1,\cdots,\mathcal{K}_{n+1},z)\\
=&(-1)^n(\theta_R(y_{n+1},z)g(\mathcal{K}_1,\cdots,\mathcal{K}_{n},x_{n+1})-\theta_R(x_{n+1},z)g(\mathcal{K}_1,\cdots,\mathcal{K}_{n},y_{n+1}))\\
&+\sum_{k=1}^{n+1}(-1)^{k+1}D_R(\mathcal{K}_k)g(\mathcal{K}_1,\cdots,\widehat{\mathcal{K}_{k}}\cdots,\mathcal{K}_{n+1},z)\\
&+\sum_{1\leq k<l\leq n+1}(-1)^k g(\mathcal{K}_1,\cdots,\widehat{\mathcal{K}_{k}}\cdots,\{x_k,y_k,x_l\}_R\wedge y_l+x_l\wedge \{x_k,y_k,y_l\}_R,\cdots,\mathcal{K}_{n+1},z)\\
&+\sum_{k=1}^{n+1}(-1)^kg(\mathcal{K}_1,\cdots,\widehat{\mathcal{K}_{k}}\cdots,\mathcal{K}_{n+1},\{x_k,y_k,z\}_R),
\end{align*}
where, $n\geq 1$,  $(f,g)\in \mathcal{C}_{\mathrm{LY}}^{n+1}(\mathfrak{L}_R,V)$, $ \mathcal{K}_i=x_i\wedge y_i\in \wedge^2 \mathfrak{L}_R, (i=1,2,\cdots,n+1), z\in \mathfrak{L}_R$.
For   any $f\in \mathcal{C}_{\mathrm{LY}}^1(\mathfrak{L}_R,V)$, its coboundary map
 $\partial^1=(\partial^1_I,\partial^1_{II})\text{:}$\, $\mathcal{C}_{\mathrm{LY}}^1(\mathfrak{L}_R,V)\rightarrow \mathcal{C}_{\mathrm{LY}}^2(\mathfrak{L}_R,V),f\rightarrow (\partial^1_I(f),\partial^1_{II}(f))$ is given by:
\begin{align*}
\partial^1_I(f)(x,y)=&\rho_R(x)f(y)-\rho_R(y)f(x)-f([x,y]_R),\\
\partial^1_{II}(f)(x,y,z)=&D_R(x,y)f(z)+\theta_R(y,z)f(x)-\theta_R(x,z)f(y)-f(\{x,y,z\}_R).
\end{align*}

\begin{definition}
Let $(V; \rho,   \theta, D, R_V)$  be a representation of the  modified  Rota-Baxter  Lie-Yamaguti algebra $(\mathfrak{L},  [-, -],   \{-, -, -\}, R)$.
Then the cochain complex $(\mathcal{C}_{\mathrm{LY}}^{\ast}(\mathfrak{L}_R,V),\partial)$ is called the cochain complex
of modified Rota-Baxter operator $R$ with coefficients in $(V; \rho_R,   \theta_R, D_R,R_V)$, denoted by $(\mathcal{C}_{\mathrm{MRBO}}^{\ast}(\mathfrak{L},V),\partial)$.
The cohomology of $(\mathcal{C}_{\mathrm{MRBO}}^{\ast}(\mathfrak{L},V),\partial)$, denoted by $\mathcal{H}_{\mathrm{MRBO}}^{\ast}(\mathfrak{L},V)$, are called the cohomology of modified
Rota-Baxter operator $R$ with coefficients in $(V; \rho_R,   \theta_R, D_R, R_V)$.
\end{definition}

In particular, when $(\mathfrak{L};\mathrm{ad}_R,\mathcal{L}_R,\mathcal{R}_R,R)$ is the adjoint  representation of $(\mathfrak{L}_R, R)$, we denote $(\mathcal{C}_{\mathrm{MRBO}}^{\ast}(\mathfrak{L},\mathfrak{L}),\partial)$ by $(\mathcal{C}_{\mathrm{MRBO}}^{\ast}(\mathfrak{L}),\partial)$
and call it the cochain complex of modified Rota-Baxter operator $R$, and denote $\mathcal{H}_{\mathrm{MRBO}}^{\ast}(\mathfrak{L},\mathfrak{L})$ by
$\mathcal{H}_{\mathrm{MRBO}}^{\ast}(\mathfrak{L})$ and call it the cohomology of modified Rota-Baxter operator $R$.

Next, we will combine the cohomology of Lie-Yamaguti algebras and the cohomology of
modified Rota-Baxter operators to construct a cohomology theory for modified Rota-Baxter Lie-Yamaguti algebras.

\begin{definition}
Let $(V; \rho,   \theta, D, R_V)$  be a representation of the  modified  Rota-Baxter  Lie-Yamaguti algebra $(\mathfrak{L},  [-, -],   \{-, -, -\}, R)$.
For any  $(f,g)\in \mathcal{C}_{\mathrm{LY}}^{n+1}(\mathfrak{L},V)$, $n\geq 1$, we define a linear map $\Phi^{n+1}=(\Phi_I^{n+1},\Phi^{n+1}_{II}):$
$\mathcal{C}^{n+1}_{\mathrm{LY}}(\mathfrak{L},V)\rightarrow \mathcal{C}^{n+1}_{\mathrm{MRBO}}(\mathfrak{L},V),  (f,g)\mapsto (\Phi_I^{n+1}(f), \Phi_{II}^{n+1}(g))$  by:
\begin{small}
\begin{align*}
&\Phi^{n+1}_{I}(f)\\
&\ =\sum_{i=1}^{n+1}\big(\sum_{1\leqslant j_1<j_2<\cdots< j_{2i-2}\leqslant {2n}}f\circ(\mathrm{Id}^{j_1-1},R,\mathrm{Id}^{j_2-j_1-1},R,\cdots, R,\mathrm{Id}^{2n-j_{2i-2}})\circ((\mathrm{Id}\wedge \mathrm{Id})^{n})\\
&\ \ \ \ -\sum_{1\leqslant j_1<j_2<\cdots< j_{2i-3}\leqslant {2n}}R_V\circ f\circ(\mathrm{Id}^{j_1-1},R,\mathrm{Id}^{j_2-j_1-1},R,\cdots, R,\mathrm{Id}^{2n-j_{2i-3}})\circ((\mathrm{Id}\wedge \mathrm{Id})^{n})\big),\\
&\Phi^{n+1}_{II}(g)\\
&\ =\sum_{i=1}^{n+1}\big(\sum_{1\leqslant j_1<j_2<\cdots< j_{2i-1}\leqslant {2n+1}}g\circ(\mathrm{Id}^{j_1-1},R,\mathrm{Id}^{j_2-j_1-1},R,\cdots, R,\mathrm{Id}^{2n+1-j_{2i-1}})\circ((\mathrm{Id}\wedge \mathrm{Id})^{n}\otimes \mathrm{Id})\\
&\ \ \ \ -\sum_{1\leqslant j_1<j_2<\cdots< j_{2i-2}\leqslant {2n+1}}R_V\circ g\circ(\mathrm{Id}^{j_1-1},R,\mathrm{Id}^{j_2-j_1-1},R,\cdots, R,\mathrm{Id}^{2n+1-j_{2i-2}})\circ((\mathrm{Id}\wedge \mathrm{Id})^{n}\otimes \mathrm{Id})\big),
\end{align*}
\end{small}
among them, when the subscript of $j_{2i-3}$ is negative, $f$ is a zero map.
In particular, when $n=0$, define $\Phi^{1}:\mathcal{C}^{1}_{\mathrm{LY}}(\mathfrak{L},V)\rightarrow \mathcal{C}^{1}_{\mathrm{MRBO}}(\mathfrak{L},V)$ by
$$\Phi^1(f)=f\circ R-R_V\circ f.$$
\end{definition}

\begin{lemma}\label{lemma:chain map}
The map $\Phi^{n+1}:\mathcal{C}^{n+1}_{\mathrm{LY}}(\mathfrak{L},V)\rightarrow \mathcal{C}^{n+1}_{\mathrm{MRBO}}(\mathfrak{L},V)$ is a cochain map, i.e., the following commutative diagram:
$$\aligned
\xymatrix{
 \mathcal{C}^1_{\mathrm{LY}}(\mathfrak{L} ,V)\ar[r]^-{\delta^1}\ar[d]^-{\Phi^1}&  \mathcal{C}^2_{\mathrm{LY}}(\mathfrak{L},V)\ar@{.}[r]\ar[d]^-{\Phi^2}& \mathcal{C}^{n+1}_{\mathrm{LY}}(\mathfrak{L},V)\ar[r]^-{\delta^{n+1}}\ar[d]^-{\Phi^{n+1}}& \mathcal{C}^{n+2}_{\mathrm{LY}}(\mathfrak{L},V)\ar[d]^{\Phi^{n+2}}\ar@{.}[r]&\\
 \mathcal{C}^1_{{\mathrm{MRBO}}}(\mathfrak{L},V)\ar[r]^-{\partial^1}& \mathcal{C}^2_{{\mathrm{MRBO}}}(\mathfrak{L},V)\ar@{.}[r]&  \mathcal{C}^{n+1}_{{\mathrm{MRBO}}}(\mathfrak{L},V)\ar[r]^-{\partial^{n+1}}& \mathcal{C}^{n+2}_{{\mathrm{MRBO}}}(\mathfrak{L},V)\ar@{.}[r]&.}
 \endaligned$$
\end{lemma}

 \begin{proof}
Because of  space limitations, here we only prove the case of~$\Phi^2\circ\delta^1=\partial^1\circ\Phi^1$.
For any $f\in \mathcal{C}^1_{\mathrm{LY}}(\mathfrak{L} ,V)$ and $x,y\in \mathfrak{L},$    we have
\noindent \begin{align}
&\Phi_I^2(\delta_I^1f)(x,y)\nonumber\\
&=(\delta_I^1f)(Rx,Ry)-R_V((\delta_I^1f)(Rx,y)+(\delta_I^1f)(x,Ry))+(\delta_I^1f)(x,y)\nonumber\\
&=\rho(Rx)f(Ry)-\rho(Ry)f(Rx)-f([Rx,Ry]))-R_V(\rho(Rx)f(y)-\rho(y)f(Rx)-f([Rx,y]))\nonumber\\
&~~~~~-R_V(\rho(x)f(Ry)-\rho(Ry)f(x)-f([x,Ry]))+\rho(x)f(y)-\rho(y)f(x)-f([x,y])\nonumber\\
&=\rho(Rx)f(Ry)-\rho(Ry)f(Rx)-f([Rx,Ry]))-R_V(\rho(Rx)f(y))+R_V(\rho(y)f(Rx))+R_V(f([Rx,y]))\nonumber\\
&~~~~~-R_V(\rho(x)f(Ry))+R_V(\rho(Ry)f(x))+R_V(f([x,Ry]))+\rho(x)f(y)-\rho(y)f(x)-f([x,y]) \label{5.1}
\end{align}
and
\noindent \begin{align}
&\partial_I^1(\Phi_I^1f)(x,y)\nonumber\\
&=\rho_R(x)(\Phi_I^1f)(y)-\rho_R(y)(\Phi_I^1f)(x)-(\Phi_I^1f)([x,y]_R)\nonumber\\
&=\rho_R(x)(f(Ry)-R_Vf(y))-\rho_R(y)(f(Rx)-R_Vf(x))-(\Phi_I^1f)([Rx,y]+[x,Ry])\nonumber\\
&=\rho(Rx)f(Ry)-R_V(\rho(x)f(Ry))-\rho(Rx)R_Vf(y)+R_V(\rho(x)R_Vf(y))-\rho(Ry)f(Rx)\nonumber\\
&~~~~~+R_V(\rho(y)f(Rx))+\rho(Ry)R_Vf(x)-R_V(\rho(y)R_Vf(x))\nonumber\\
&~~~~~-f(R[Rx,y])+R_V(f[Rx,y])-f(R[x,Ry])+R_V(f[x,Ry]).\label{5.2}
\end{align}
Using   Eqs. \eqref{3.1}, \eqref{4.1}, \eqref{4.6} and \eqref{4.8}, further comparing  Eqs. \eqref{5.1} and \eqref{5.2}, we have \eqref{5.1}=\eqref{5.2}.
Therefore, $\Phi_I^2\circ\delta_I^1=\partial_I^1\circ\Phi_I^1$. Similarly, there is also  $\Phi_{II}^2\circ\delta_{II}^1=\partial_{II}^1\circ\Phi_{II}^1$.
\end{proof}

\begin{definition}
Let $(V; \rho,   \theta, D, R_V)$  be a representation of the  modified  Rota-Baxter  Lie-Yamaguti algebra $(\mathfrak{L},  [-, -],   \{-, -, -\}, R)$.
We define the set of 1-cochains of  modified  Rota-Baxter  Lie-Yamaguti algebras by  $\mathcal{C}^{1}_{\mathrm{MRBLY}}(\mathfrak{L},V)=\mathcal{C}^{1}_{\mathrm{LY}}(\mathfrak{L},V)$.
  For $n\geq 1$, we define the set of  $(n+1)$-cochains of modified  Rota-Baxter  Lie-Yamaguti algebras by
$$\mathcal{C}_{\mathrm{MRBLY}}^{n+1}(\mathfrak{L},V):=
\mathcal{C}^{n+1}_{\mathrm{LY}}(\mathfrak{L},V)\oplus \mathcal{C}^{n}_{\mathrm{MRBO}}(\mathfrak{L},V).$$

Define  a linear map  $\mathrm{d}^1:\mathcal{C}_{\mathrm{MRBLY}}^{1}(\mathfrak{L},V)\rightarrow \mathcal{C}_{\mathrm{MRBLY}}^{2}(\mathfrak{L},V)$  by
\begin{align*}
\mathrm{d}^1(f_1)=(\delta^1 f_1, -\Phi^1 (f_1)), \forall f_1\in \mathcal{C}_{\mathrm{MRBLY}}^{1}(\mathfrak{L},V).
\end{align*}
For  $n = 1$,   we define  the linear map  $\mathrm{d}^{2}:\mathcal{C}_{\mathrm{MRBLY}}^{2}(\mathfrak{L},V)\rightarrow \mathcal{C}_{\mathrm{MRBLY}}^{3}(\mathfrak{L},V)$  by
\begin{align*}
\mathrm{d}^{2}((f_1,g_1),f_2)=(\delta^{2}(f_1,g_1), -\partial^{1} (f_2)- \Phi^{2}(f_1,g_1)),
\end{align*}
for any  $((f_1,g_1),f_2)\in \mathcal{C}_{\mathrm{MRBLY}}^{2}(\mathfrak{L},V).$\\
Then, for  $n \geq 2$,   we define  the linear map  $\mathrm{d}^{n+1}:\mathcal{C}_{\mathrm{MRBLY}}^{n+1}(\mathfrak{L},V)\rightarrow \mathcal{C}_{\mathrm{MRBLY}}^{n+2}(\mathfrak{L},V)$  by
\begin{align*}
\mathrm{d}^{n+1}((f_1,g_1),(f_2,g_2))=(\delta^{n+1}(f_1,g_1), -\partial^{n} (f_2,g_2)- \Phi^{n+1}(f_1,g_1)),
\end{align*}
for any  $((f_1,g_1),(f_2,g_2))\in \mathcal{C}_{\mathrm{MRBLY}}^{n+1}(\mathfrak{L},V).$
\end{definition}

In view of Lemma \ref{lemma:chain map}, we have the following theorem.

\begin{theorem}\label{theorem: cochain complex}
The map $\mathrm{d}^{n+1}$ is a coboundary operator, i.e., $\mathrm{d}^{n+1}\circ\mathrm{d}^{n}=0.$
\end{theorem}

Therefore, from Theorem \ref{theorem: cochain complex}, we obtain a cochain complex $(\mathcal{C}_{\mathrm{MRBLY}}^{\ast}(\mathfrak{L},V),\mathrm{d})$   called the cochain complex
of modified  Rota-Baxter  Lie-Yamaguti algebra $(\mathfrak{L},  [-, -],   \{-, -, -\}, R)$ with coefficients in $(V; \rho,   \theta, D, R_V)$.
The cohomology of  $(\mathcal{C}_{\mathrm{MRBLY}}^{\ast}(\mathfrak{L},V),\mathrm{d})$, denoted by
 $\mathcal{H}_{\mathrm{MRBLY}}^{\ast}(\mathfrak{L},V)$,
 is called the cohomology of the  modified  Rota-Baxter  Lie-Yamaguti algebra  $(\mathfrak{L},  [-, -],   \{-, -, -\}, R)$ with coefficients in $(V; \rho,   \theta, D, R_V)$.

 In particular,  when $(V; \rho,   \theta, D, R_V)=(\mathfrak{L};\mathrm{ad},\mathcal{L},\mathcal{R},R)$, we just denote $(\mathcal{C}_{\mathrm{MRBLY}}^{\ast}(\mathfrak{L},\mathfrak{L}),\mathrm{d})$, $\mathcal{H}_{\mathrm{MRBLY}}^{\ast}(\mathfrak{L},\mathfrak{L})$  by
 $(\mathcal{C}_{\mathrm{MRBLY}}^{\ast}(\mathfrak{L}),\mathrm{d})$, $\mathcal{H}_{\mathrm{MRBLY}}^{\ast}(\mathfrak{L})$ respectively, and call them the cochain complex, the cohomology of modified  Rota-Baxter  Lie-Yamaguti algebra $(\mathfrak{L},  [-, -], \{-, $  $ -, -\}, R)$ respectively.

It is obvious that there is a  short exact sequence of cochain complexes:
\begin{align*}
0\rightarrow \mathcal{C}_{\mathrm{MRBO}}^{\ast-1}(\mathfrak{L},V)\stackrel{}{\longrightarrow}\mathcal{C}_{\mathrm{MRBLY}}^{\ast}(\mathfrak{L},V)\stackrel{}{\longrightarrow}\mathcal{C}_{\mathrm{LY}}^{\ast}(\mathfrak{L},V)\rightarrow 0.
\end{align*}
It induces a long exact sequence of cohomology groups:
\begin{align*}
\cdots\rightarrow \mathcal{H}_{\mathrm{MRBLY}}^{p}(\mathfrak{L},V)\rightarrow \mathcal{H}_{\mathrm{LY}}^{p}(\mathfrak{L},V)\rightarrow \mathcal{H}_{\mathrm{MRBO}}^{p}(\mathfrak{L},V)\rightarrow \mathcal{H}_{\mathrm{MRBLY}}^{p+1}(\mathfrak{L},V)\rightarrow \mathcal{H}_{\mathrm{LY}}^{p+1}(\mathfrak{L},V)\rightarrow \cdots.
\end{align*}

 \section{  Formal deformations of modified Rota-Baxter Lie-Yamaguti algebras}\label{sec:deformations}
\def\theequation{\arabic{section}. \arabic{equation}}
\setcounter{equation} {0}

Motivated by the  deformations of Lie-Yamaguti algebras\cite{L15,Z15}, in this section, we study   formal deformations of modified Rota-Baxter Lie-Yamaguti algebras.
Let $\mathbb{K}[[t]]$ be a ring of power series of one variable $t$, and let $\mathfrak{L}[[t]]$ be the set of formal power series over $\mathfrak{L}$.  If
$(\mathfrak{L}, [-, -], \{-, -, -\})$ is a Lie-Yamaguti algebra, then there is a Lie-Yamaguti algebra structure over the
ring $\mathbb{K}[[t]]$ on $\mathfrak{L}[[t]]$ given by
\begin{align*}
&[\sum_{i=0}^{\infty}x_it^i,\sum_{i=0}^{\infty}y_jt^j]=\sum_{s=0}^{\infty}\sum_{i+j=s}[x_i,y_j]t^s,\\
& \{\sum_{i=0}^{\infty}x_it^i,\sum_{i=0}^{\infty}y_jt^j,\sum_{k=0}^{\infty}z_kt^k\}=\sum_{s=0}^{\infty}\sum_{i+j+k=s}\{x_i,y_j,z_k\}t^s.
\end{align*}


\begin{definition}
A   formal deformation of the modified Rota-Baxter  Lie-Yamaguti algebra  $(\mathfrak{L}, [-, -], \{-, -, -\},R)$ is a triple  $ (F_t, G_t, R_t)$  of the forms
$$F_t=\sum_{i=0}^{\infty}F_it^i,~~ G_t=\sum_{i=1}^{\infty}G_it^i,~~R_t=\sum_{i=0}^{\infty}R_it^i,$$
such that the following conditions are satisfied:

(i) $((F_i,G_i),R_i)\in \mathcal{C}^{2}_{\mathrm{MRBLY}}(\mathfrak{L});$

(ii) $F_0=[-, -],G_0=\{-, -, -\}$ and $R_0=R;$

(iii) and $(\mathfrak{L}[[t]], F_t, G_t, R_t)$  is  a   modified Rota-Baxter  Lie-Yamaguti algebra over $\mathbb{K}[[t]]$.
\end{definition}

Let $ (F_t, G_t, R_t)$  be a formal deformation as above. Then, for any $ x, y, z, a, b\in \mathfrak{L}$, the following equations must hold:
\begin{align*}
&~F_t(x, y)+F_t(y,  x)=0,~~G_t(x,y,z)+G_t(y,x,z)=0,\\
&~ \circlearrowleft_{x,y,z}F_t(F_t(x, y), z)+\circlearrowleft_{x,y,z}G_t(x,y,z)=0,\\
&~\circlearrowleft_{x,y,z}G_t(F_t(x, y), z, a)=0,\\
&~G_t(a, b, F_t(x, y))=F_t(G_t(a, b, x), y)+F_t(x,G_t(a, b, y)),\\
& ~ G_t(a, b, G_t(x, y, z))=G_t(G_t(a, b, x), y, z)+ G_t(x,  G_t(a, b, y), z)+ G_t(x,  y, G_t(a, b, z)),\\
&~ F_t(R_tx, R_ty)=R_t(F_t(R_tx,y)+F_t(x,R_ty))-F_t(x,y), \\
&~G_t(R_tx,R_ty,R_tz)=R_t(G_t(x,R_ty,R_tz)+G_t(R_tx,y,R_tz)+G_t(R_tx,R_ty,z)+G_t(x,y,z)) \\
&~~~~~~~~~~~~~~~~~~~~~~~~~~~~~-G_t(R_tx,y,z)-G_t(x,R_ty,z)-G_t(x,y,R_tz).
\end{align*}
Collecting the coefficients of $t^n$, we get that the above equations are equivalent to the following equations.
\begin{align}
&~F_n(x, y)+F_n(y,  x)=0,~~G_n(x,y,z)+G_n(y,x,z)=0,\label{6.1}\\
& \sum_{i=0}^n\circlearrowleft_{x,y,z}F_i(F_{n-i}(x, y), z)+\circlearrowleft_{x,y,z}G_n(x,y,z)=0,\label{6.2}\\
&\sum_{i=0}^n\circlearrowleft_{x,y,z}G_i(F_{n-i}(x, y), z, a)=0,\label{6.3}\\
&\sum_{i=0}^nG_i(a, b, F_{n-i}(x, y))=\sum_{i=0}^nF_i(G_{n-i}(a, b, x), y)+\sum_{i=0}^nF_i(x,G_{n-i}(a, b, y)),\label{6.4}\\
& \sum_{i=0}^n G_i(a, b, G_{n-i}(x, y, z))\nonumber\\
&=\sum_{i=0}^n\big(G_i(G_{n-i}(a, b, x), y, z)+ G_i(x,  G_{n-i}(a, b, y), z)+ G_i(x,  y, G_{n-i}(a, b, z))\big),\label{6.5}\\
&~ \sum_{i+j+k=n}F_i(R_jx, R_ky)=\sum_{i+j+k=n}R_i(F_j(R_kx,y)+F_j(x,R_ky))-F_n(x,y), \label{6.6}\\
&\sum_{i+j+k+l=n}G_i(R_jx,R_ky,R_lz)\nonumber\\
\quad &=\sum_{i+j+k+l=n}R_i\big(G_j(x,R_ky,R_lz)+G_j(R_kx,y,R_lz)+G_j(R_kx,R_ly,z)\big)\nonumber\\
\quad &+\sum_{i+j=n}\big(R_i(G_j(x,y,z))-G_i(R_jx,y,z)-G_i(x,R_jy,z)-G_i(x,y,R_jz)\big).\label{6.7}
\end{align}
Note that for $n=0$, equations \eqref{6.1}--\eqref{6.7}  are equivalent to $(\mathfrak{L}, F_0, G_0, R_0)$  is  a   modified Rota-Baxter  Lie-Yamaguti algebra.

\begin{proposition}\label{prop:2-cocycle}
Let $ (F_t, G_t, R_t)$  be a formal deformation of  a modified Rota-Baxter  Lie-Yamaguti algebra  $(\mathfrak{L}, [-, -], \{-, -, -\},R)$.
Then $((F_1,G_1),R_1)$ is a 2-cocycle in the cochain complex $(\mathcal{C}^{\ast}_{\mathrm{MRBLY}}(\mathfrak{L}),\mathrm{d})$.
\end{proposition}
\begin{proof}
For  $n=1$,   equations \eqref{6.1}--\eqref{6.7} become
\begin{align}
&F_1(x, y)+F_1(y,  x)=0,~~G_1(x,y,z)+G_1(y,x,z)=0,\label{6.8}\\
&[F_{1}(x, y), z]+F_1([x, y], z)+[F_{1}(z, x), y]+F_1([z, x], y)+[F_{1}(y, z), x]+F_1([y, z], x)\nonumber\\
&+G_1(x,y,z)+G_1(z,x,y)+G_1(y,z, x)=0,\label{6.9}\\
&G_1([x, y], z, a)+\{F_{1}(x, y), z, a\}+G_1([z, x], y, a)+\{F_{1}(z, x), y, a\}+G_1([y, z], x, a)\nonumber\\
&+\{F_{1}(y, z), x, a\}=0,\label{6.10}\\
&G_1(a, b, [x, y])+\{a, b, F_{1}(x, y)\}\nonumber\\
&=F_1(\{a, b, x\}, y)+[G_{1}(a, b, x), y]+F_1(x,\{a, b, y\})+[x,G_{1}(a, b, y)],\label{6.11}\\
&  G_1(a, b, \{x, y, z\})+  \{a, b, G_{1}(x, y, z)\}=G_1(\{a, b, x\}, y, z)+ \{G_{1}(a, b, x), y, z\}\nonumber\\
&+ G_1(x,  \{a, b, y\}, z)+ \{x,  G_{1}(a, b, y), z\}+ G_1(x,  y, \{a, b, z\})+ \{x,  y, G_{1}(a, b, z)\},\label{6.12}\\
&F_1(Rx, Ry)+[R_1x, Ry]+[Rx, R_1y]=R_1([Rx,y]+[x,Ry])+R(F_1(Rx,y)+F_1(x,Ry))\nonumber\\
&+R([R_1x,y]+[x,R_1y])-F_1(x,y), \label{6.13}\\
&G_1(Rx,Ry,Rz)+\{R_1x,Ry,Rz\}+\{Rx,R_1y,Rz\}+\{Rx,Ry,R_1z\}\nonumber\\
&=R_1\big(\{x,Ry,Rz\}+\{Rx,y,Rz\}+\{Rx,Ry,z\}\big)+R\big(G_1(x,Ry,Rz)+G_1(Rx,y,Rz)\nonumber\\
&+G_1(Rx,Ry,z)\big)+R\big(\{x,R_1y,Rz\}+\{R_1x,y,Rz\}+\{R_1x,Ry,z\}\nonumber\\
&+\{x,Ry,R_1z\}+\{Rx,y,R_1z\}+\{Rx,R_1y,z\}\big)\nonumber\\
& +R_1(\{x,y,z\})+R(G_1(x,y,z))-G_1(Rx,y,z)-G_1(x,Ry,z)-G_1(x,y,Rz)\nonumber\\
&-\{R_1x,y,z\}-\{x,R_1y,z\}-\{x,y,R_1z\}.\label{6.14}
\end{align}
From Eqs. \eqref{6.8}--\eqref{6.12}, we get $(F_1,G_1)\in \mathcal{C}^{2}_{\mathrm{LY}}(\mathfrak{L})$ and $\delta^{2}(F_1,G_1)=0$.  Further from    Eqs. \eqref{6.13} and \eqref{6.14}, we have
$-\partial_I^{1} (R_1)- \Phi_I^{2}(F_1)=0$ and $-\partial_{II}^{1} (R_1)- \Phi_{II}^{2}(G_1)=0$  respectively.
 Hence, $\mathrm{d}^{2}((F_1,G_1),R_1)=0,$ that is, $((F_1,G_1),R_1)$ is a 2-cocycle in   $(\mathcal{C}^{\ast}_{\mathrm{MRBLY}}(\mathfrak{L}),\mathrm{d})$.
\end{proof}

\begin{definition}
The 2-cocycle  $((F_1,G_1),R_1)$ is called the infinitesimal of the   formal deformation $ (F_t, G_t, R_t)$ of a modified Rota-Baxter  Lie-Yamaguti algebra  $(\mathfrak{L}, [-, -], \{-, -, -\},R)$.
\end{definition}

\begin{definition}
Let $(F_t, G_t, R_t)$ and $(F'_t, G'_t, R'_t)$ be two formal deformations of
a  modified Rota-Baxter  Lie-Yamaguti algebra  $(\mathfrak{L}, [-, -], \{-, -, -\},R)$. A formal isomorphism from $(\mathfrak{L}[[t]], F_t, G_t, R_t)$ to $(\mathfrak{L}[[t]],F'_t, G'_t, R'_t)$ is a power series
$\Psi_t=\sum_{i=o}^{\infty}\Psi_it^i:\mathfrak{L}[[t]]\rightarrow \mathfrak{L}[[t]]$ , where  $\Psi_i:\mathfrak{L}\rightarrow \mathfrak{L}$ are linear maps
with $\Psi_0=\mathrm{id}_\mathfrak{L}$, such that:
\begin{align}
& \Psi_t \circ F_t=F'_t\circ(\Psi_t\otimes\Psi_t),\label{6.15}\\
& \Psi_t \circ G_t=G'_t\circ(\Psi_t\otimes\Psi_t\otimes\Psi_t),\label{6.16}\\
& \Psi_t \circ R_t=R'_t\circ \Psi_t,\label{6.17}
\end{align}
In this case, we say that the two   formal deformations $(F_t, G_t, R_t)$ and $(F'_t, G'_t, R'_t)$
are equivalent.
\end{definition}

\begin{proposition}
The infinitesimals of two equivalent   formal deformations of $(\mathfrak{L}, [-, -],$ $ \{-, -, -\},R)$
are in the same cohomology class in $\mathcal{H}^{2}_{\mathrm{MRBLY}}(\mathfrak{L})$.
\end{proposition}
\begin{proof}
Let $\Psi_t:(\mathfrak{L}[[t]], F_t, G_t, R_t)\rightarrow (\mathfrak{L}[[t]],F'_t, G'_t, R'_t)$ be a formal isomorphism.
 By expanding Eqs. \eqref{6.15}--\eqref{6.17} and comparing the coefficients of $t$ on both sides, we have
 \begin{align*}
&F'_1-F_1=F_0\circ (\Psi_1\o \mathrm{Id}_\mathfrak{L})+F_0\circ (\mathrm{Id}_\mathfrak{L} \o\Psi_1)-\Psi_1\circ F_0,\\
&G'_1-G_1=G_0\circ (\Psi_1\o \mathrm{Id}_\mathfrak{L}\o \mathrm{Id}_\mathfrak{L})+G_0\circ (\mathrm{Id}_\mathfrak{L}\o  \mathrm{Id}_\mathfrak{L} \o\Psi_1)+G_0\circ (\mathrm{Id}_\mathfrak{L}\o \Psi_1 \o \mathrm{Id}_\mathfrak{L})-\Psi_1\circ G_0,\\
&R'_{1}-R_1= R \circ \Psi_1- \Psi_{1}\circ R.
 \end{align*}
That is, we get
  $$((F'_1,G'_1), R'_1)-((F_1,G_1), R_1)=(\delta^1(\Psi_1),-\Phi^1(\Psi_1))=\mathrm{d}^1(\Psi_1)\in \mathcal{B}_{\mathrm{MRBLY}}^{2}(\mathfrak{L}). $$
  Therefore,  $((F'_1,G'_1), R'_1)$ and $((F_1,G_1), R_1)$ are in the same cohomology class in $\mathcal{H}^{2}_{\mathrm{MRBLY}}(\mathfrak{L})$.
\end{proof}

\begin{definition}
Given a modified Rota-Baxter  Lie-Yamaguti algebra  $(\mathfrak{L}, [-, -], \{-, -, -\},R)$, the power series $F_t, G_t, R_t$ with $F_i=\mathrm{I}_{i,0}F_0, G_i=\mathrm{I}_{i,0}G_0, R_i=\mathrm{I}_{i,0}R_0$ make
$(F_t, G_t, R_t)$  into a   formal deformation of $(\mathfrak{L}, [-, -], \{-, -, -\},R)$. A formal
deformation  equivalent to this one is called trivial  deformation.
\end{definition}

\begin{definition}
A modified Rota-Baxter  Lie-Yamaguti algebra  $(\mathfrak{L}, [-, -], \{-, -, -\},R)$ is said to be rigid if every   formal deformation of $\mathfrak{L}$ is trivial deformation.
\end{definition}

\begin{theorem}
 Let $(\mathfrak{L}, [-, -], \{-, -, -\},R)$ be a modified Rota-Baxter  Lie-Yamaguti algebra. If $\mathcal{H}^{2}_{\mathrm{MRBLY}}(\mathfrak{L})=0$, then
$(\mathfrak{L}, [-, -], \{-, -, -\},R)$ is rigid.
\end{theorem}
\begin{proof}
Let $(F_t, G_t, R_t)$ be a   formal deformation of $(\mathfrak{L}, [-, -], \{-, -, -\},R)$. From Proposition \ref{prop:2-cocycle},
$((F_1, G_1), R_1)$ is a 2-cocycle. By $\mathcal{H}^{2}_{\mathrm{MRBLY}}(\mathfrak{L})=0$, there exists a 1-cochain
$$\Psi_1 \in \mathcal{C}^1_{\mathrm{MRBLY}}(\mathfrak{L})=\mathcal{C}^1_{\mathrm{LY}}(\mathfrak{L})$$
such that $((F_1,G_1), R_1)=  \mathrm{d}^1(\Psi_1), $  that is, $F_1=\delta_I^1(\Psi_1)$, $G_1=\delta_{II}^1(\Psi_1)$ and $R_1=-\Phi^1(\Psi_1)$.

Setting $\Psi_t = \mathrm{id}_\mathfrak{L} -\Psi_1t$, we get a deformation  $(\overline{F}_t, \overline{G}_t, \overline{R}_t)$, where
 \begin{align*}
\overline{F}_t=&\Psi_t^{-1}\circ F_t\circ (\Psi_t\otimes \Psi_t),\\
\overline{G}_t=&\Psi_t^{-1}\circ G_t\circ (\Psi_t\otimes \Psi_t\otimes \Psi_t),\\
\overline{R}_t=&\Psi_t^{-1}\circ R_t\circ \Psi_t.
 \end{align*}
 It can be easily verify that $\overline{F}_1=0,\overline{G}_1=0,  \overline{R}_1=0$. Then
    $$\begin{array}{rcl} \overline{F}_t&=& F_0+\overline{F}_2t^2+\cdots,\\
    \overline{G}_t&=& G_0+\overline{G}_2t^2+\cdots,\\
 R_t&=& R_0+\overline{R}_2t^2+\cdots.\end{array}$$
  By Eqs.     \eqref{6.1}--\eqref{6.7},
  we see that $((\overline{F}_2, \overline{G}_2),  \overline{R}_2)$  is still a 2-cocyle, so by induction, we can show
that $(F_t, G_t, R_t)$  is equivalent to the trivial deformation $(F_0, G_0, R_0)$.  Therefore,  $(\mathfrak{L}, [-, -], \{-, -, -\},R)$ is rigid.
\end{proof}

  \section{   Abelian extensions of modified Rota-Baxter Lie-Yamaguti algebras}\label{sec:extensions}
\def\theequation{\arabic{section}. \arabic{equation}}
\setcounter{equation} {0}

In this section, we consider  abelian extensions of modified Rota-Baxter Lie-Yamaguti algebras.
We prove that any abelian extension of a modified Rota-Baxter Lie-Yamaguti algebra  has a representation and a  2-cocycle.  It is further proved that they are classified by the second cohomology.

Notice that a vector space $V$ together with a linear map $R_V:V\rightarrow V$ is naturally an
abelian modified Rota-Baxter Lie-Yamaguti algebra where the multiplications on  $V$ is defined to be $[-,-]_V=0,\{-,-,-\}_V=0.$

\begin{definition}
 Let $(\mathfrak{L}, [-, -], \{-, -, -\},R)$  be a modified Rota-Baxter Lie-Yamaguti algebra and $(V, [-, -]_V, \{-, -, -\}_V,R_V)$  an abelian
modified Rota-Baxter Lie-Yamaguti algebra.
An abelian extension of  $(\mathfrak{L}, [-, -], \{-, -, -\},R)$ by  $(V, [-, -]_V, \{-, -, -\}_V,R_V)$
 is  a short exact sequence of   morphisms of modified Rota-Baxter Lie-Yamaguti algebras
$$\begin{tiny}\begin{CD}
0@>>> {(V, [-, -]_V, \{-, -, -\}_V,R_V)} @>i >> (\hat{\mathfrak{L}},[-, -]_{\hat{\mathfrak{L}}}, \{-, -, -\}_{\hat{\mathfrak{L}}},\hat{R}) @>p >> (\mathfrak{L}, [-, -], \{-, -, -\},R) @>>>0
\end{CD}\end{tiny}$$
such that  $V$ is an abelian ideal of $\hat{\mathfrak{L}},$  i.e.,  $[u, v]_{\hat{\mathfrak{L}}}=0,\{-, u,v\}_{\hat{\mathfrak{L}}}=\{u,v,-\}_{\hat{\mathfrak{L}}}=0, \forall u,v\in V$.
\end{definition}

\begin{definition}
 A   section  of an abelian extension $(\hat{\mathfrak{L}},[-, -]_{\hat{\mathfrak{L}}}, \{-, -, -\}_{\hat{\mathfrak{L}}},\hat{R})$ of $(\mathfrak{L}, [-, -],$ $ \{-, -, -\},R)$  by  $(V, [-, -]_V, \{-, -, -\}_V,R_V)$ is a linear map $s:\mathfrak{L}\rightarrow \hat{\mathfrak{L}}$ such that   $p\circ s=\mathrm{id}_\mathfrak{L}$.
\end{definition}

\begin{definition}
   Let $(\hat{\mathfrak{L}}_1,[-, -]_{\hat{\mathfrak{L}}_1}, \{-, -, -\}_{\hat{\mathfrak{L}}_1},\hat{R}_1)$ and  $(\hat{\mathfrak{L}}_2,[-, -]_{\hat{\mathfrak{L}}_2}, \{-, -, -\}_{\hat{\mathfrak{L}}_2},\hat{R}_2)$  be two abelian extensions of $(\mathfrak{L}, [-, -],$ $ \{-, -, -\},R)$  by  $(V, [-, -]_V, \{-, -, -\}_V,R_V)$. They are said to be  equivalent if  there is an isomorphism of  modified Rota-Baxter Lie-Yamaguti algebrass $\varphi:(\hat{\mathfrak{L}}_1,[-, -]_{\hat{\mathfrak{L}}_1}, \{-, -, -\}_{\hat{\mathfrak{L}}_1},\hat{R}_1)\rightarrow (\hat{\mathfrak{L}}_2,[-, -]_{\hat{\mathfrak{L}}_2}, \{-, -, -\}_{\hat{\mathfrak{L}}_2},\hat{R}_2)$
such that the following diagram is  commutative:
\begin{tiny}
\begin{align}
\begin{CD}
0@>>> {(V, [-, -]_V, \{-, -, -\}_V,R_V)} @>i_1 >> (\hat{\mathfrak{L}}_1,[-, -]_{\hat{\mathfrak{L}}_1}, \{-, -, -\}_{\hat{\mathfrak{L}}_1},\hat{R}_1) @>p_1 >> (\mathfrak{L}, [-, -],  \{-, -, -\},R) @>>>0\\
@. @| @V \varphi VV @| @.\\
0@>>> {(V, [-, -]_V, \{-, -, -\}_V,R_V)} @>i_2 >> (\hat{\mathfrak{L}}_2,[-, -]_{\hat{\mathfrak{L}}_2}, \{-, -, -\}_{\hat{\mathfrak{L}}_2},\hat{R}_2) @>p_2 >> (\mathfrak{L}, [-, -],  \{-, -, -\},R) @>>>0.\label{7.1}
\end{CD}
\end{align}
\end{tiny}
\end{definition}

Now for an  abelian extension $(\hat{\mathfrak{L}},[-, -]_{\hat{\mathfrak{L}}}, \{-, -, -\}_{\hat{\mathfrak{L}}},\hat{R})$ of $(\mathfrak{L}, [-, -],$ $ \{-, -, -\},R)$  by  $(V, [-, -]_V, \{-, -, -\}_V,R_V)$  with a section $s:\mathfrak{L}\rightarrow\hat{\mathfrak{L}}$, we define  linear maps $\rho: \mathfrak{L} \rightarrow \mathfrak{gl}(V)$  and  $\theta, D: \mathfrak{L}\times\mathfrak{L} \rightarrow \mathfrak{gl}(V)$  by
$$\varrho(x)u:=[s(x),u]_{\hat{\mathfrak{L}}},  $$
$$\theta(x,y)u:=\{u,s(x),s(y)\}_{\hat{\mathfrak{L}}},  \quad \forall x,y\in \mathfrak{L}, u\in V.$$
In particular, $D(x,y)u=\theta(y,x)u-\theta(x,y)u=\{s(x),s(y),u\}_{\hat{\mathfrak{L}}}.$
\begin{proposition} \label{prop:representation}
  With the above notations, $(V; \rho, \theta,D, R_V)$ is a representation of the  modified Rota-Baxter Lie-Yamaguti algebra  $(\mathfrak{L}, [-, -],$ $ \{-, -, -\},R)$.
\end{proposition}
\begin{proof}
 In view of \cite{Z15}, $(V; \rho, \theta,D)$ is a representation of a  Lie-Yamaguti algebra  $(\mathfrak{L}, [-, -],$ $ \{-, -, -\})$.
 Further, for any $x,y\in \mathfrak{L}$ and $ u\in V,$ ${\hat{R}}s(x)-s(Rx)\in V$ means that $\rho({\hat{R}}s(x))u=\rho(s(Rx))u,\theta({\hat{R}}s(x), {\hat{R}}s(y)) u=\theta(s(Rx), s(Ry)) u$. Therefore,
we have
 \begin{align*}
&\rho(Rx)R_Vu=[s(Rx),R_Vu]_{\hat{\mathfrak{L}}}=[\hat{R}s(x),R_Vu]_{\hat{\mathfrak{L}}}\\
&=R_V([\hat{R}s(x),u]_{\hat{\mathfrak{L}}}+[s(x),R_Vu]_{\hat{\mathfrak{L}}})-[s(x),u]_{\hat{\mathfrak{L}}}\\
&=R_V(\rho(Rx)u+\rho(x)R_Vu)-\rho(x)u,\\
&\theta(Rx,  Ry)R_Vu\\
\quad  &=  [R_V u, s(Rx), s(Ry)]_{\hat{\mathfrak{L}}}\\
\quad  &= [R_V u, {\hat{R}}s(x), {\hat{R}}s(y)]_{\hat{\mathfrak{L}}}\\
\quad  &=R_V( \{R_V u, {\hat{R}}s(x), s(y)\}_{\hat{\mathfrak{L}}}+\{u, {\hat{R}}s(x), {\hat{R}}s(y)\}_{\hat{\mathfrak{L}}}+\{R_V u, s(x), {\hat{R}}s(y)\}_{\hat{\mathfrak{L}}}+\{u, s(x), s(y)\}_{\hat{\mathfrak{L}}})\\
\quad  &~~~~~~- \{R_V u, s(x), s(y)\}_{\hat{\mathfrak{L}}}-\{u, {\hat{R}}s(x), s(y)\}_{\hat{\mathfrak{L}}}-\{u, s(x), {\hat{R}}s(y)\}_{\hat{\mathfrak{L}}}\\
\quad  &=R_V( \{R_V u, s(Rx), s(y)\}_{\hat{\mathfrak{L}}}+\{u, s(Rx), s(Ry)\}_{\hat{\mathfrak{L}}}+\{R_V u, s(x), s(Ry)\}_{\hat{\mathfrak{L}}}+\{u, s(x), s(y)\}_{\hat{\mathfrak{L}}})\\
\quad  &~~~~~~- \{R_V u, s(x), s(y)\}_{\hat{\mathfrak{L}}}-\{u, s(Rx), s(y)\}_{\hat{\mathfrak{L}}}-\{u, s(x), s(Ry)\}_{\hat{\mathfrak{L}}}\\
\quad  &=R_V\big(\theta(Rx, Ry)u +\theta(Rx, y)R_Vu+\theta(x, Ry)R_Vu+\theta(x, y)u\big)\\
\quad &~~~~~~~- \theta(Rx, y)u-\theta(x, y)R_Vu-\theta(x, Ry)u.
\end{align*}
Hence,  $(V; \rho, \theta,D, R_V)$ is a representation of   $(\mathfrak{L}, [-, -],$ $ \{-, -, -\},R)$.
\end{proof}

We further define linear maps $\nu:\mathfrak{L}\times \mathfrak{L}\rightarrow V$, $\psi:\mathfrak{L}\times \mathfrak{L}\times \mathfrak{L}\rightarrow V$ and $\chi:\mathfrak{L}\rightarrow V$ respectively by
\begin{align*}
\nu(x,y)&=[s(x), s(y)]_{\hat{\mathfrak{L}}}-s([x, y]),\\
\psi(x,y,z)&=\{s(x), s(y), s(z)\}_{\hat{\mathfrak{L}}}-s(\{x, y,z\}),\\
\chi(a)&=\hat{R}(s(x))-s(R(x)),\quad\forall x,y,z\in \mathfrak{L}.
\end{align*}
We transfer the modified Rota-Baxter Lie-Yamaguti algebra structure on $\hat{\mathfrak{L}}$ to $\mathfrak{L}\oplus V$ by endowing $\mathfrak{L}\oplus V$ with   multiplications $[-, -]_\nu,\{-,-,-\}_\psi$
and a modified Rota-Baxter
operator  $R_\chi$  defined by
\begin{align}
[x+u, y+v]_\nu&=[x, y]+\rho(x)v-\rho(y)u+\nu(x,y), \label{7.2}\\
 \{x+u, y+v, z+w\}_\psi&=\{x, y, z\}+\theta(y, z)u-\theta(x, z)v+D(x, y)w+\psi(x, y, z),\label{7.3}\\
 R_\chi(x+u)&=R(x)+\chi(x)+R_V(u),  \forall x,y,z\in \mathfrak{L},\,u,v,w\in V.\label{7.4}
\end{align}

\begin{proposition}\label{prop:2-cocycles}
The 4-tuple $(\mathfrak{L}\oplus V,[-,-]_\nu, \{-,-,-\}_\psi, R_\chi)$ is a modified Rota-Baxter Lie-Yamaguti algebra  if and only if
$((\nu,\psi),\chi)$ is a 2-cocycle  of the modified Rota-Baxter Lie-Yamaguti algebra $(\mathfrak{L}, [-, -], \{-, -, -\},R)$ with the coefficient  in $(V; \rho, \theta,D, R_V)$.
 In this case,
$$ \begin{tiny}\begin{CD}
0@>>> {(V, [-, -]_V, \{-, -, -\}_V,R_V)} @>i >> (\mathfrak{L}\oplus V,[-,-]_\nu, \{-,-,-\}_\psi, R_\chi) @>p >>(\mathfrak{L}, [-, -], \{-, -, -\},R) @>>>0
\end{CD}\end{tiny}$$
 is an abelian extension.
\end{proposition}
\begin{proof}
In view of \cite{Z15},
$(\mathfrak{L}\oplus V,[-,-]_\nu, \{-,-,-\}_\psi)$ is a  Lie-Yamaguti algebra  if and only if
$\delta^2(\nu,\psi)=0$.

The map $R_\chi$ is a modified modified Rota-Baxter  operator on $(\mathfrak{L}\oplus V,[-,-]_\nu, \{-,-,-\}_\psi)$ if and only~if
\begin{align}
&[R_\chi(x + u), R_\chi(y + v)]_\nu\nonumber\\
\quad &=R_\chi([R_\chi(x + u), y + v]_\nu + [x + u, R_\chi(y + v)]_\nu+[x + u, y + v]_\nu)-[x + u, y + v]_\nu,\label{7.5}\\
&\{R_\chi(x + u), R_\chi(y + v), R_\chi(z + w)\}_\psi\nonumber\\
\quad &=R_\chi(\{R_\chi(x + u), R_\chi(y + v), z + w\}_\psi + \{x + u, R_\chi(y + v), R_\chi(z + w)\}_\psi\nonumber\\
\quad & +\{R_\chi(x + u), y + v,R_\chi(z + w)\}_\psi)+\{x + u, y + v, z + w\}_\psi)-\{R_\chi(x + u), y + v, z + w\}_\psi\nonumber\\
\quad &-\{x + u, R_\chi(y + v), z + w\}_\psi - \{x + u, y + v, R_\chi(z + w)\}_\psi,\label{7.6}
\end{align}
for any $x,y,z\in \mathfrak{L}$ and $u,v,w\in V$. Further,
Eqs  \eqref{7.5} and \eqref{7.6} are equivalent to the following equations:
\begin{align}
&\nu(Rx, Ry)+\rho(Rx)\chi(y)-\rho(Ry)\chi(x)=\chi([Rx,y])+\chi([x,Ry])+R_V(\rho(x)\chi(y))\nonumber\\
&-R_V(\rho(y)\chi(x))+R_V(\nu(Rx,y))+R_V(\nu(x,Ry))-\nu(x,y) \label{7.7}\\
&\psi(Rx, Ry, Rz) + \theta(Ry, Rz)\chi(x)-\theta(Rx, Rz)\chi(y) + D(Rx, Ry)\chi(z)\nonumber\\
\quad &=R_V (\psi(Rx, Ry, z) + \psi(Rx, y, Rz) + \psi(x, Ry, Rz) + \psi(x, y, z)) -\psi(Rx, y, z) - \psi(x, Ry, z)\nonumber\\
\quad &-\psi(x, y, Rz) - \theta(Ry, z)\chi(x)+ \theta(Rx, z)\chi(y)-\theta(y, Rz)\chi(x)-D(Rx, y)\chi(z)\nonumber\\
 \quad &+\theta(x, Rz)\chi(y) -D(x, Ry)\chi(z)+R_V(\theta(y, z)\chi(x)-\theta(x, z)\chi(y)+D(x, y)\chi(z))\nonumber\\
\quad & + \chi(\{Rx, Ry, z\} + \{Rx, y, Rz\} + \{x, Ry, Rz\}+\{x, y, z\}).\label{7.8}
\end{align}
Using Eqs.  \eqref{7.7} and  \eqref{7.8}, we get $-\partial_I^{1} (\chi)- \Phi_I^{2}(\nu)=0$ and $-\partial_{II}^{1} (\chi)- \Phi_{II}^{2}(\psi)=0$  respectively.
Therefore, $ \mathrm{d}^2((\nu,\psi),\chi)=(\delta^2(\nu,\psi),-\partial^1(\chi)-\Phi^{2}(\nu,\psi))=0,$ that is, $((\nu,\psi),\chi)$ is a  2-cocycle.

Conversely, if $((\nu,\psi),\chi)$ is a  2-cocycle  of the modified Rota-Baxter Lie-Yamaguti algebra $(\mathfrak{L}, [-, -], \{-, -, -\},R)$ with the coefficient  in $(V; \rho, \theta,D, R_V)$,  then we have $ \mathrm{d}^2((\nu,\psi),\chi)=(\delta^2(\nu,\psi),-\partial^1(\chi)-\Phi^{2}(\nu,\psi))=0,$ in which $\delta^2(\nu,\psi)=0$, Eqs.  \eqref{7.7} and \eqref{7.8} hold.
So $(\mathfrak{L}\oplus V,[-,-]_\nu, \{-,-,-\}_\psi, R_\chi)$  is a  modified Rota-Baxter Lie-Yamaguti algebra.
\end{proof}

\begin{proposition}
Let $(\hat{\mathfrak{L}},[-, -]_{\hat{\mathfrak{L}}}, \{-, -, -\}_{\hat{\mathfrak{L}}},\hat{R})$ be an abelian extension of $(\mathfrak{L}, [-, -],$ $ \{-, -, -\},R)$  by  $(V, [-, -]_V, \{-, -, -\}_V,R_V)$ and $s:\mathfrak{L}\rightarrow\hat{\mathfrak{L}}$   a section. If   $((\nu,\psi),\chi)$ is a 2-cocycle   constructed using the section $s$, then  its cohomology class does not depend on the choice of $s$.
\end{proposition}
\begin{proof}
Let $s_1,s_2:\mathfrak{L}\rightarrow \hat{ \mathfrak{L}}$ be two distinct sections, then we have two corresponding 2-cocycles $((\nu_1,\psi_1),\chi_1)$ and $((\nu_2,\psi_2),\chi_2)$ respectively. Define
a linear map $\lambda: \mathfrak{L}\rightarrow V$ by $\lambda(x)=s_1(x)-s_2(x)$. Then
\begin{align*}
&\nu_1(x, y) \\
\quad &= [s_1(x), s_1(y)]_{\hat{\mathfrak{L}}}-s_1[x, y]\\
\quad &= [s_2(x) + \lambda(x), s_2(y) + \lambda(y)]_{\hat{\mathfrak{L}}}- s_2([x, y])-\lambda[x, y]\\
\quad &= [s_2(x), s_2(y)]_{\hat{\mathfrak{L}}} +\rho(x)\lambda(y)-\rho(y)\lambda(x)-s_2([x, y])- \lambda[x, y]\\
\quad &= \nu_2(x, y) + \delta_{I}^1\lambda(x, y),\\
&\psi_1(x, y, z) \\
\quad &= \{s_1(x), s_1(y), s_1(z)\}_{\hat{\mathfrak{L}}}-s_1\{ x, y, z\}\\
\quad &= \{s_2(x) + \lambda(x), s_2(y) + \lambda(y), s_2(z) + \lambda(z)\}_{\hat{\mathfrak{L}}}- s_2(\{x, y, z\})-\lambda\{x, y, z\})\\
\quad &= \{s_2(x), s_2(y), s_2(z)\}_{\hat{\mathfrak{L}}} + \theta(y, z)\lambda(x)-\theta(x, z)\lambda(y)\\
\quad &+D(x, y)\lambda(z)-s_2(\{x, y, z\})- \lambda(\{x, y, z\})\\
\quad &= \psi_2(x, y, z) + \delta_{II}^1\lambda(x, y, z),
\end{align*}
and
\begin{align*}
\chi_1(x) &= \hat{R}s_1(x)- s_1R(x)\\
\quad &= \hat{R}(s_2(x) + \lambda(x))-(s_2(R(x)) + \lambda(Rx)\\
\quad &= \hat{R}s_2(x)-s_2R(x) + \hat{R}\lambda(x)-\lambda(Rx)\\
\quad &= \chi_2(x) + R_V \lambda(x)-\lambda (Rx)\\
\quad &= \chi_2(x)-\Phi^1\lambda(x).
\end{align*}
So, $((\nu_1,\psi_1),\chi_1)-((\nu_2,\psi_2),\chi_2)=(\delta^1\lambda,-\Phi^1\lambda)=\mathrm{d}^1(\lambda)\in  \mathcal{B}_{\mathrm{MRBLY}}^{2}(\mathfrak{L},V)$, that is  $((\nu_1,\psi_1),\chi_1)$ and $((\nu_2,\psi_2),\chi_2)$ are in the same cohomological class  in $ \mathcal{H}_{\mathrm{MRBLY}}^{2}(\mathfrak{L},V)$.
\end{proof}

\begin{theorem} \label{theorem:classify abelian extensions}
Abelian extensions of a modified Rota-Baxter Lie-Yamaguti algebra  $(\mathfrak{L}, [-, -],$ $ \{-, -, -\},R)$  by  $(V, [-, -]_V, \{-, -, -\}_V,R_V)$ are classified by the second cohomology group $\mathcal{H}_{\mathrm{MRBLY}}^{2}(\mathfrak{L},V)$ of $(\mathfrak{L}, [-, -],$ $ \{-, -, -\},R)$ with coefficients in the representation $(V; \rho, \theta, D, R_V)$.
\end{theorem}
\begin{proof}
Suppose $(\hat{\mathfrak{L}}_1,[-, -]_{\hat{\mathfrak{L}}_1}, \{-, -, -\}_{\hat{\mathfrak{L}}_1},\hat{R}_1)$ and  $(\hat{\mathfrak{L}}_2,[-, -]_{\hat{\mathfrak{L}}_2}, \{-, -, -\}_{\hat{\mathfrak{L}}_2},\hat{R}_2)$
 are equivalent abelian extensions   of $(\mathfrak{L}, [-, -],$ $ \{-, -, -\},R)$  by  $(V, [-, -]_V, \{-, -, -\}_V,R_V)$  with the associated isomorphism $\varphi:(\hat{\mathfrak{L}}_1,[-, -]_{\hat{\mathfrak{L}}_1}, \{-, -, -\}_{\hat{\mathfrak{L}}_1},\hat{R}_1)\rightarrow (\hat{\mathfrak{L}}_2,[-, -]_{\hat{\mathfrak{L}}_2}, \{-, -, -\}_{\hat{\mathfrak{L}}_2},\hat{R}_2)$ such that the diagram in~\eqref{7.1} is commutative.
 Let $s_1$ be a section of $(\hat{\mathfrak{L}}_1,[-, -]_{\hat{\mathfrak{L}}_1}, \{-, -, -\}_{\hat{\mathfrak{L}}_1},\hat{R}_1)$. As $p_2\circ f=p_1$, we  have
$$p_2\circ(\varphi\circ s_1)=p_1\circ s_1= \mathrm{id}_{\mathfrak{L}}.$$
That is, $\varphi\circ s_1$ is a section of $(\hat{\mathfrak{L}}_2,[-, -]_{\hat{\mathfrak{L}}_2}, \{-, -, -\}_{\hat{\mathfrak{L}}_2},\hat{R}_2)$. Denote $s_2:=\varphi\circ s_1$. Since $\varphi$ is an isomorphism of  modified Rota-Baxter Lie-Yamaguti algebra such that $\varphi|_V=\mathrm{id}_V$, we have
\begin{align*}
\nu_2(x,y)&=[s_2(x), s_2(y)]_{\hat{\mathfrak{L}}_2}-s_2([x,y])\\
&=[\varphi(s_1(x)), \varphi(s_1(y))]_{\hat{\mathfrak{L}}_2}-\varphi(\sigma_1([x, y]))\\
&=\varphi\big([s_1(x), s_1(y)]_{\hat{\mathfrak{L}}_1}-s_1([x, y])\big)\\
&=\varphi(\nu_1(x,y))\\
&=\nu_1(x,y),\\
\psi_2(x, y, z)&=\{s_2(x), s_2(y), s_2(z)\}_{\hat{\mathfrak{L}}_2}-s_2(\{x, y, z\})\\
\quad &= \varphi(\{s_1(x), s_1(y), s_1(z)\}_{\hat{\mathfrak{L}}_1}-s_1\{x, y, z\})\\
\quad &= \psi_1(x, y, z),\\
\chi_2(x)&=\hat{R}(s_2(x))-s_2(R(x))\\
\quad &=\hat{R}(\varphi\circ s_1(x))-\varphi\circ s_1(R(x))\\
\quad &= \hat{R}(s_1(x))- s_1(R(x))\\
\quad &=\chi_1(x).
\end{align*}
So, all equivalent abelian extensions give rise to the same element in $\mathcal{H}_{\mathrm{MRBLY}}^{2}(\mathfrak{L},V)$.

Conversely, given two  cohomologous 2-cocycles $((\nu_1,\psi_1),\chi_1)$ and $((\nu_2,\psi_2),\chi_2)$  in  $\mathcal{H}_{\mathrm{MRBLY}}^{2}(\mathfrak{L},V)$,
we can construct two abelian extensions $(\mathfrak{L}\oplus V,[-,-]_{\nu_1}, \{-,-,-\}_{\psi_1}, R_{\chi_1})$ and  $(\mathfrak{L}\oplus V,[-,-]_{\nu_2}, \{-,-,-\}_{\psi_2}, R_{\chi_2})$ via Proposition \ref{prop:2-cocycles}. Then  there is  a linear map $\lambda: \mathfrak{L}\rightarrow  V$ such that
 $$((\nu_1,\psi_1),\chi_1)-((\nu_2,\psi_2),\chi_2)=\mathrm{d}^1(\lambda)=(\delta^1\lambda,-\Phi^1\lambda).$$

 Define a linear map $\varphi_\lambda: \mathfrak{L}\oplus V\rightarrow  \mathfrak{L}\oplus V$ by
$\varphi_\lambda(x+u):=x+\lambda(x)+u, ~x\in \mathfrak{L}, u\in V.$  Then we have $\varphi_\lambda$ is an isomorphism of these two abelian extensions.
\end{proof}
\begin{center}
 {\bf ACKNOWLEDGEMENT}
 \end{center}

 The paper is supported by the NSF of China (No. 12161013) and   Guizhou Provincial Basic Research Program (Natural Science) (No. ZK[2023]025).

\renewcommand{\refname}{REFERENCES}

\end{document}